\documentstyle[11pt]{article}

\newtheorem{theorem}{Theorem}

\newtheorem{proposition}[theorem]{Proposition}
\newtheorem{definition}[theorem]{Definition}
\newtheorem{example}[theorem]{Example}
\newtheorem{remark}[theorem]{Remark}

\newcommand{\defin}[1]
  {\begin{definition} {\rm #1} \end{definition}}
\newcommand{\examp}[1]
  {\begin{example} {\rm #1} \end{example}}

\def\QED{\quad\blackslug\lower 8.5pt\null}

  \newcommand{\htimes}{%
    \mathbin{\mathsurround0pt \mathchoice
      {\hbox{\vrule\negthinspace$\times$}}%
      {\hbox{\vrule\negthinspace$\times$}}%
      {\hbox{\vrule\negthinspace$\scriptstyle\times$}}%
      {\hbox{\vrule\negthinspace$\scriptscriptstyle\times$}}%
    }}

\begin{document}

\begin{center}
{\Large \bf  ON THE THEORY OF ALMOST}
 
\vspace*{2mm}

{\Large \bf  GRASSMANN STRUCTURES}

\vspace*{3mm}
{\large M.A. Akivis and  V.V. Goldberg}

\end{center}

{\footnotesize 
\noindent
{\bf Abstract}. 
The differential geometry of almost Grassmann structures defined 
on a differentiable manifold of dimension $n = pq$ by a fibration 
of Segre cones $SC (p, q)$ is studied. The peculiarities in 
the structure of  almost Grassmann 
structures for the cases $p=q=2; \;
p = 2, q > 2$ (or $p > 2, q = 2$), and $p > 2, q > 2$ 
are clarified. The fundamental geometric objects of these 
structures up to fourth order are derived. The conditions 
under which an almost Grassmann structure is locally flat or 
locally semiflat are found for all cases indicated above.
}

\setcounter{equation}{0}


\setcounter{equation}{0}

\setcounter{section}{-1}

\section{Introduction}  

Almost Grassmann manifolds were introduced by Hangan \cite{H66} 
as a generalization of the Grassmannian $G(m, n)$.  Hangan 
(\cite{H66}, \cite{H68}) and Ishihara \cite{I72} studied mostly some special almost Grassmann manifolds, especially locally 
Grassmann manifolds.  Later almost Grassmann manifolds were 
studied by Goldberg \cite{G75}, Mikhailov \cite{M78} and 
Akivis \cite{A82a} in connection with the development of 
the theory of multidimensional webs. 
Goncharov \cite{Go87} considered the almost Grassmann manifolds  
as generalized conformal structures. 

~
Baston in \cite{B91} constructed a theory of a general class 
of 
structures, called almost Hermitian symmetric (AHS) structures, 
which include conformal, projective, almost Grassmann, and 
quaternionic structures and for which the construction of the 
Cartan normal connection is possible. He constructed a tensor 
invariant for them and proved that its vanishing is equivalent to 
the structure being locally that of a Hermitian symmetric space.  
In \cite{Go87}, the AHS structures have been studied from 
the point of view of cone structures (see \cite{B91} and 
\cite{Go87} for further references on generalized conformal 
structures). 

Bailey and Eastwood \cite{BE91} extended the theory of local 
twistors, which was 
known for four-dimensional conformal structures, to the  
almost Grassmann structures (they  called them the paraconformal 
structures).   Dhooghe (\cite{D93}, \cite{D94}) 
considered almost Grassmann structures (he 
called them Grassmannian structures) as subbundles of the second 
order frame bundle and constructed a canonical normal connection 
for these structures. The structure equations derived for the 
Grassmann structures in \cite{D94} 
are very close to the structure equations of the  Grassmann 
structures considered in the present paper.

In the current paper we consider  the real theory 
of  almost Grassmann structures while in \cite{Go87}, 
\cite{B91}, and \cite{BE91} their complex theory was studied.

Although some of the authors who studied almost Grassmann  
structures proved that an almost Grassmann  structure is 
a $G$-structure of finite type two (see \cite{H80} and 
\cite{M78}), no one of the authors went further than the 
development of the  first structure tensor. 

In Section 1 of the current paper we 
consider the Grassmann structures. In Section 2 we 
 define the  almost Grassmann  
structures geometrically following \cite{A80} and \cite{A82a} 
(see also \cite{M78}, \cite{AS92}, \S 8.3, and \cite{G88}, \S\S 2.1 and 2.2). 

Sections 3 and 4 are the core of the paper. In these sections 
 we make the most symmetric and natural specializations of 
second-,  third-, and fourth-order moving frames, derive the 
structure equations and construct (in a fourth-order 
differential neighborhood)   a complete  geometric object of
 the almost Grassmann  structure totally defining its geometric 
structure, find the expression 
 of the components of the torsion tensor of an almost 
Grassmann structure in the general (not reduced) third-order 
frame, prove that an almost Grassmann structure is a 
$G$-structure of finite type two (Theorem 4.1), establish 
relations between the components of the complete geometric object 
(Theorems 4.2 and 4.3) and use these connections to determine 
under what conditions an almost Grassmann structure is 
locally Grassmann (Theorem 4.4).

In Section 4 we also find the structure group of the  almost 
Grassmann  structure these structures and its  
differential prolongation. The structure group of the  
almost Grassmann structure is represented in the form  
${\bf SL} (p) \times {\bf SL} (q) \times {\bf H}$, 
where    ${\bf SL} (p), {\bf SL} (q)$ are the special 
linear groups of order $p$ and $q$; respectively, and ${\bf H}$ 
is the group of homotheties.  The prolonged 
group $G'$  is isomorphic to the semidirect product 
$G \htimes {\bf T} (pq)$, where ${\bf T} (pq)$ is the 
 group  of parallel translations of $(pq)$-dimensional affine 
space. 

The almost Grassmann  structure defines on the manifold $M$ two 
fiber bundles $E_\alpha$ and $E_\beta$ that are real twistor 
fiber bundles. In a more general setting they were considered in 
[BE 91]. For the general almost Grassmann structure, 
the first nonvanishing structure tensor 
splits into two subtensors, which are the structure tensors of 
fiber bundles $E_\alpha$ and $E_\beta$. The vanishing of each of 
these subtensors leads to integrability of the corresponding 
twistor fiber bundle. In Section 5 we give without proof the 
statements of our main theorems on semiintegrability of almost 
Grassmann structures. Their complete proofs are given in 
our book [AG 93].

Note also the main differences between our paper and the 
papers \cite{D93} and \cite{D94}: 
\begin{enumerate}
\item 
We consider a projective model of a Grassmann manifold while in 
the papers \cite{D93} and \cite{D94} (and in most of other papers 
on this subject) a vector model was considered. This is the 
reason that in our presentation both the dimension  $n$ of the 
original projective space and the dimension $m$ of the generating 
element of the Grassmann manifold $G (m, n)$ are one unit less 
than in those papers. The dimension 
of the Grassmann manifold $G (m, n)$ and the dimension 
of the manifold $M$ endowed with an almost Grassmann  structure 
are equal to $pq$ where $p = m + 1$ and $ q = n -m$.

\item 
In this paper an almost Grassmann  structure on a manifold 
$M$ is defined geometrically by means of a fibration of 
Segre cones $SC (p, q)$ while in  \cite{D93} and \cite{D94} 
it is defined analytically by means of structure equations.

\item We clearly distinguish three cases: 
\begin{description}
\item[a)] $p = 2, q = 2 \;(m = 1, n = 3), \dim M = 4$. 
In this case the almost Grassmann  structure 
$AG (1, 3)$ is equivalent to the pseudoconformal structure 
$CO (2, 2)$ that is torsion-free and whose conformal curvature 
is determined by its Weyl tensor defined in a third-order 
differential neighborhood.

\item[b)] $p = 2, q > 2$ (or $p > 2, q = 2$). 
In this case the torsion tensor of the fiber bundle 
$E_\alpha$ (respectively, $E_\beta$) vanishes, and 
the difference between this structure and a locally flat 
 structure is determined by the curvature tensor of 
the fiber bundle $E_\alpha$ and the torsion tensor 
of $E_\beta$ (respectively, $E_\beta$ and $E_\alpha$).

\item[c)] $p > 2, q > 2$. In this case the difference between 
this structure and a locally flat  structure is completely 
determined by the torsion tensors of the fiber bundles $E_\alpha$ 
$E_\beta$.
\end{description}

\item We construct the fundamental geometric objects of these 
structures up to fourth order for each of these three cases and 
establish connections among them.

\item In the main parts of \cite{D93} and \cite{D94} the author 
considered torsion-free almost Grassmann  structures. He    
did not have results of our Theorems {\bf 4.2, 4.3} 
and {\bf 4.4} and erroneously assumed that for 
$p > 2, q > 2$ there exist torsion-free almost Grassmann  
structures that are not locally Grassmann  (locally flat) 
structures. According to our Theorem {\bf 4.4}, 
if $p > 2$ and $q > 2$, then 
torsion-free almost Grassmann  structures are locally Grassmann  
structures. This is the reason that the results 
of  the main parts of \cite{D93} and \cite{D94} valid 
only for $p = 2, q > 2; p > 2, q = 2$, and $p = q = 2$.   
\end{enumerate}

\section{Grassmann Structures}

\setcounter{equation}{0}
 
{\bf 1.} 
Let $P^n$ be an $n$-dimensional projective space.  The set of 
$m$-dimensional subspaces $P^m \subset P^n$ is called the {\em 
Grassmann manifold}, or the {\em Grassmannian}, and is denoted by 
the symbol $G (m, n)$. It is well-known that the Grassmannian is 
a differentiable manifold, and that its dimension is equal to 
$\rho = (m + 1)(n - m)$. It will be convenient for us to set 
$p = m + 1$ and $q = n - m$. Then we have $n = p + q - 1$.
 
Let a subspace $P^m = x$ be an element of the Grassmannian 
$G (m, n)$. With any subspace  $x$, we associate a family  of 
projective point frames $\{A_u\}, u = 0, 1, \ldots n$, such that 
the vertices $A_\alpha, \alpha = 0, 1, \ldots , m$, of its frames 
lie in the subspace $P^m$, and the points $A_i, i = m+1, \ldots, 
n$, lie outside $P^m$ and together with the points $A_\alpha$ make up the frame $\{A_u\}$ of the space  $P^n$. 

We will write the equations of infinitesimal displacement of the 
moving frames we have chosen in the form:
\begin{equation}\label{eq:1.1}
d A_u = \theta_u^v A_v, \;\; u, v = 0, \ldots, n.
\end{equation}
Since the fundamental group  of the space $P^n$ is locally 
isomorphic to the group ${\bf SL} (n+1)$, the forms 
$\theta_u^v$ are connected by the relation 
\begin{equation}\label{eq:1.2}
\theta_u^u = 0. 
\end{equation}
The structure equations of the space $P^n$ have the form
\begin{equation}\label{eq:1.3}
d \theta_u^v = \theta_u^w \wedge \theta_w^v. 
\end{equation}
By (1.3), the exterior differential of the left-hand side of 
equations (1.1) is identically equal to 0, and hence the system 
of equations (1.1) as well as equation (1.2) is completely integrable.

By (1.1), we have 
$$
d A_\alpha = \theta_\alpha^\beta A_\beta + \theta_\alpha^i A_i.
$$
It follows that the 1-forms $\theta_\alpha^i$ are basis forms of 
the Grassmannian. These forms are linearly independent, and their 
number is equal to $\rho = (m+1)(n-m) = p \cdot q$, i.e., it 
equals the dimension of the 
Grassmannian $G (m, n)$. We will assume  that the integers $p$ 
and $q$ satisfy the inequalities $p \geq 2$ and $q \geq 2$, since 
for $p = 1$, we have $m = 0$, and the Grassmannian $G (0, n)$ is 
the projective space $P^n$, and for $q=1$, we have $m = n - 1$, 
and the Grassmannian $G (n-1, n)$ is isomorphic to the dual 
projective space $(P^n)^*$.

Let us rename the basis forms by setting $\theta_\alpha^i = 
\omega_\alpha^i$ and finding their exterior differentials:
\begin{equation}\label{eq:1.4}
d \omega_\alpha^i = \theta_\alpha^\beta \wedge \omega_\beta^i 
+ \omega_\alpha^j \wedge \theta^i_j. 
\end{equation}
Define the trace-free forms 
\begin{equation}\label{eq:1.5}
\omega_\alpha^\beta = \theta_\alpha^\beta 
- \frac{1}{p} \delta_\alpha^\beta \theta^\gamma_\gamma, \;\;
\omega^i_j = \theta^i_j  
- \frac{1}{q} \delta_j^i  \theta^k_k,
\end{equation}
satisfying the conditions:
\begin{equation}\label{eq:1.6}
\omega_\alpha^\alpha = 0, \;\;\; \omega_i^i = 0.
\end{equation}
 Eliminating the forms $\theta_\alpha^\beta $ and $\theta_i^j$ 
from equations (1.4), we find that
\begin{equation}\label{eq:1.7}
d \omega_\alpha^i = \omega_\alpha^\beta \wedge \omega_\beta^i 
+ \omega_\alpha^j \wedge \omega_j^i + \omega \wedge \omega_\alpha^i, 
\end{equation}
where $\omega = \displaystyle \frac{1}{p} \theta^\gamma_\gamma 
- \displaystyle \frac{1}{q} \theta^k_k$, or by (1.2), 
\begin{equation}\label{eq:1.8}
\omega = \Bigl(\displaystyle \frac{1}{p} + \displaystyle  
\frac{1}{q}\Bigr)   \theta^\gamma_\gamma.
\end{equation}

Setting
\begin{equation}\label{eq:1.9}
 \omega_i^\alpha = - \Bigl(\displaystyle \frac{1}{p} + \displaystyle  \frac{1}{q}\Bigr) \theta_i^\alpha    
\end{equation}
and taking the exterior derivatives of equations (1.7) and (1.8), 
we obtain 
\begin{equation}\label{eq:1.10}
\renewcommand{\arraystretch}{1.3}
\begin{array}{ll}
d \omega_\alpha^\beta = \omega_\alpha^\gamma \wedge 
 \omega_\gamma^\beta + \displaystyle \frac{q}{p+q} 
\omega_\gamma^k \wedge 
(\delta_\alpha^\beta \omega_k^\gamma - p \delta_\alpha^\gamma 
\omega_k^\beta), \\
d \omega_j^i = \omega_j^k \wedge \omega_k^i + 
\displaystyle \frac{p}{p+q} (\delta^i_j \omega_k^\gamma 
- q \delta^i_k \omega_j^\gamma) 
\wedge \omega_\gamma^k  
\end{array}
\renewcommand{\arraystretch}{1}
\end{equation}
and 
\begin{equation}\label{eq:1.11}
 d \omega = \omega_i^\alpha \wedge \omega^i_\alpha.    
\end{equation}
Exterior differentiation of equations (1.9) gives 
\begin{equation}\label{eq:1.12}
 d \omega_i^\alpha  = \omega_i^j \wedge \omega_j^\alpha 
+ \omega_i^\beta \wedge \omega_\beta^\alpha + 
\omega_i^\alpha \wedge \omega.    
\end{equation}
Finally,  exterior differentiation of equations (1.12) leads to 
identities.

Thus, the structure equations of the Grassmannian $G (m, n)$ take 
the form (1.7), (1.10),  (1.11) and (1.12). This system of 
differential equations is {\em closed} in the sense that its 
further exterior differentiation leads to identities. 

If we fix a subspace $x = P^m \subset P^n$, then we obtain 
$\omega^i_\alpha = 0$, and equations (1.10) and (1.11) become
\begin{equation}\label{eq:1.13}
d \pi_\alpha^\beta = \pi_\alpha^\gamma \wedge 
 \pi_\gamma^\beta, \;\; d \pi_j^i = \pi_j^k \wedge \pi_k^i, \;\;  
d \pi = 0,  
\end{equation}
where  $ \pi = \omega (\delta), 
\pi_\alpha^\beta = \omega_\alpha^\beta (\delta), 
\pi_j^i = \omega_j^i (\delta)$, 
and $\delta$ is the operator of differentiation with respect to 
the fiber parameters of the second-order frame bundle associated 
with the Grassmannian $G (m, n)$. Moreover, the forms 
$\pi_\alpha^\beta$ and $\pi_j^i$ satisfy  equations similar to 
equations (1.6), that is, these forms are trace-free. The forms 
$\pi_\alpha^\beta$ are invariant forms of the group 
${\bf SL} (p)$ which is locally isomorphic to the group of 
projective transformations of the subspace $P^m$. The forms 
$\pi_j^i$ are invariant forms of the group ${\bf SL} (q)$ which 
is locally isomorphic to the group of projective transformations 
of the bundle  of $(m+1)$-dimensional subspaces of the space 
$P^n$ containing $P^m$. The form $\pi$ is an invariant form of 
the group ${\bf H} = {\bf R}^* \otimes \mbox{{\rm Id}}$ of 
homotheties of the 
space $P^n$ with center at $P^m$; here ${\bf R}^*$ is the 
multiplicative group of real numbers. 

The direct product of these three groups is the structural group 
$G$ of the Grassmann manifold $G (m, n)$:
\begin{equation}\label{eq:1.14}
G = {\bf SL} (p) \times {\bf SL} (q) \times {\bf H}.
\end{equation}

Finally, the forms $\pi_i^\alpha = \omega_i^\alpha (\delta)$, 
which by (1.12) satisfy the structure equations 
\begin{equation}\label{eq:1.15}
d \pi_i^\alpha = \pi_i^j \wedge  \pi_j^\alpha + \pi_i^\beta 
\wedge \pi_\beta^\alpha + \pi_i^\alpha \wedge \pi,   
\end{equation}
are also fiber forms on the Grassmannian $G (m, n)$ but unlike 
the forms $\pi_\alpha^\beta,\; \pi_i^j$ and $\pi$, they are 
connected with the third-order frame bundle of the Grassmannian 
$G (m, n)$. 

The forms  $\pi_\alpha^\beta, \;\pi_j^i, \pi$ and 
$\pi_i^\alpha$, satisfying the structure equations 
(1.13) and (1.15), are invariant forms of the  group 
\begin{equation}\label{eq:1.16}
G' = G \htimes {\bf T} (pq)
\end{equation}
arising under the differential prolongation 
of the structure group $G$ of the Grassmannian $G (m, n)$. 
The group $G'$ is the group of motions of 
an $(n - m - 1)$-quasiaffine space 
$A^n_{n-m-1}$ (see \cite{Do88}) which is 
a projective space $P^n$ with a fixed  $m$-dimensional subspace 
$P^m = A_0 \wedge A_1 \wedge \ldots \wedge A_m$ and 
the generating element 
$P^{n-m-1} = A_{m+1} \wedge \ldots \wedge A_n$. 
The dimension of the space 
$A^n_{n-m-1}$ coincides with the 
dimension  of the Grassmannian $G (n - m - 1, n)$, and this 
dimension is the same as the dimension of the Grassmannian 
$G (m, n)$: $\rho = (m + 1)(n - m)$. 
The forms  $\pi_i^\alpha$ are invariant forms of the  group  
${\bf T} (pq)$  of parallel translations of 
the space $A^n_{n-m-1}$, and  the group $G$  is the 
stationary subgroup  of its element $P^{n-m-1}$.

In the index-free notation, the structure equations (1.7) and 
(1.10)--(1.12) of the Grassmannian $G (m, n)$ can be written as 
follows:
\begin{equation}\label{eq:1.17}
\renewcommand{\arraystretch}{1.3}
\begin{array}{ll}
 d\theta =  \omega \wedge \theta -  \theta  \wedge \omega_\alpha 
-  \omega_\beta \wedge \theta, \\
d \omega_\alpha  + \omega_\alpha \wedge \omega_\alpha  
= \displaystyle \frac{q}{p+q} \Bigl[- I_\alpha 
\mbox{{\rm tr}}\; (\varphi \wedge \theta) 
+ p \varphi \wedge \theta\Bigr],  \\
d \omega_\beta  + \omega_\beta \wedge \omega_\beta  
= \displaystyle \frac{p}{p+q} \Bigl[- I_\beta 
\mbox{{\rm tr}}\; (\varphi \wedge \theta) 
+ q \theta \wedge \varphi\Bigr],  \\
d \omega  =  \mbox{{\rm tr}} \; (\varphi \wedge \theta), \\
 d \varphi  + \theta_\alpha \wedge \varphi  + \varphi \wedge 
\theta_\beta  + \omega \wedge  \varphi = 0,
\end{array}
\renewcommand{\arraystretch}{1}
\end{equation}
where $\theta = (\omega_\alpha^i)$ is the matrix 1-forms 
defined in the first-order fiber bundle;  
$\omega_\alpha = (\omega_\beta^\alpha)$ and 
$\omega_\beta = (\omega_j^i)$ are the matrix 1-forms 
defined in a second-order fiber bundle for which 
$$
\mbox{{\rm tr}} \; \omega_\alpha = 0, \;\; 
\mbox{{ \rm tr}} \; \omega_\beta = 0;
$$
 the form $\omega$ is the scalar form  occurring in equations 
(1.7), (1.11), and (1.12) 
and also defined in a second-order frame bundle;  
$\varphi = (\omega_i^\alpha)$ is a matrix 1-form defined in a 
third-order fiber bundle; and 
 $I_\alpha = (\delta_\alpha^\beta)$ and 
$I_\beta = (\delta_i^j)$ are the unit tensors of orders $p$ and 
$q$, respectively. 

 Note that   in the exterior products 
of 1-forms, occurring in equations (1.17) and in 
further structure equations of this subsection, 
  multiplication is performed  according 
to the regular rules of matrix multiplication---row by  column.

Along with the Grassmannian $G (m, n)$, in the space $P^n$ 
one can consider the dual manifold $G (n - m - 1, n)$. Its base 
forms are the forms $\omega_i^\alpha$, and its geometry is identical to that of the Grassmannian $G (m, n)$.

{\bf 2.} With the help of Grassmann coordinates, the Grassmannian 
$G (m, n)$ can be mapped onto a smooth algebraic variety 
$\Omega (m, n)$ of dimension $\rho = pq$ embedded into a 
projective space $P^N$ of dimension 
$N = \displaystyle {p+q \choose p} - 1$.

Suppose that $x = A_0 \wedge A_1 \wedge \ldots \wedge A_m $ 
is a point of the variety $\Omega (m, n)$. Then 
\begin{equation}\label{eq:1.18}
d x = \tau x + \omega_\alpha^i e^\alpha_i,
\end{equation}
where by (1.8) $\tau = \theta_0^0 + \ldots + \theta_m^m 
= \frac{pq}{p+q} \omega$ and 
$$
e_i^\alpha = A_0 \wedge \ldots  \wedge A_{\alpha-1} 
\wedge A_i  \wedge A_{\alpha+1}  \wedge \ldots  \wedge A_m, 
$$
and the points $e_i^\alpha$ together with the point $x$ determine 
a basis in the tangent subspace $T_x (\Omega)$. The second 
differential of the point $x$ satisfies the relation
\begin{equation}\label{eq:1.19}
d^2 x \equiv \sum_{\alpha<\beta, i< j} 
(\omega_{\alpha}^{i} \omega_{\beta}^{j} - \omega_{\alpha}^{j} \omega_{\beta}^{i})  e^{\alpha \beta}_{ij} \pmod{T_x (\Omega)},
\end{equation}
where  
$$
 e^{\alpha \beta}_{ij} =  A_0 \wedge \ldots  \wedge A_{\alpha-1} 
 \wedge A_i  \wedge A_{\alpha+1}  \wedge \ldots \wedge 
A_{\beta-1} \wedge A_j  \wedge A_{\beta+1}  \wedge \ldots \wedge 
A_m 
$$
are points of the space $P^N$ that lie  on the variety 
$\Omega (m, n)$ and together with the points $x$ and $e_i^\alpha$ 
determine the osculating space $T_x^2 (\Omega)$ 
of the variety $\Omega$. The quadratic forms 
\begin{equation}\label{eq:1.20}
\omega_{\alpha\beta}^{ij} = \omega_{\alpha}^{i} 
\omega_{\beta}^{j} - \omega_{\alpha}^{j} \omega_{\beta}^{i}  
\end{equation}
are the second fundamental forms of the variety $\Omega \subset 
P^N$. 

The equations $\omega_{\alpha\beta}^{ij} = 0$ determine the cone 
of asymptotic directions of the variety $\Omega$ at a point 
$x \in \Omega$. The equations of this cone can be written as 
follows:
\begin{equation}\label{eq:1.21}
{\rm rank}\;\; (\omega^i_\alpha) = 1.
\end{equation}
In view of (1.21), parametric equations of this cone have 
the form 
\begin{equation}\label{eq:1.22}
\omega_{\alpha}^{i} = t_{\alpha} s^{i}, \;\; 
\alpha = 0, 1, \ldots , m;\;\; i = m+1, \ldots , n. 
\end{equation}

If we consider a projectivization of this cone, then $t_\alpha$ 
and $s^i$ can be taken as homogeneous coordinates of projective 
spaces $P^m$ and $P^{n-m-1}$. Thus such a projectivization is an 
embedding of the direct product $P^{p-1} \times P^{q-1}$ into a 
projective space $P^{\rho - 1}$ of dimension ${\rho - 1}$, 
where $\rho = pq$. Such 
an embedding is called  the {\em Segre variety} and is denoted by $S (p-1, q-1)$. This is the reason that the cone of asymptotic 
directions of the variety $\Omega$ determined by equations (1.21) 
is called the {\em Segre cone}. This cone is denoted by 
$SC_x (p, q)$ since it carries two families of plane generators 
of dimensions $p$ and $q$.  Plane generators from different 
families of the cone $SC_x (p, q)$ have a common straight line. 
It is possible to prove that the cone $SC_x (p, q)$ is the 
intersection of the tangent subspace $T_x (\Omega)$ and the 
variety $\Omega$:
$$
C_x (p, q) = T_x (\Omega) \cap \Omega.
$$

The differential geometry of Grassmannians was studied  
in the papers \cite{A82b}, \cite{L61}, and \cite{W67}.

\section{Almost Grassmann Structures}

\setcounter{equation}{0}
 
{\bf 1.} 
Now we can define the notion of an almost Grassmann structure. 

\defin{\label{def:2.1} Let $M$ be a differentiable manifold of 
dimension $p q$, and let $SC (p, q)$ be a differentiable 
fibration of Segre cones with the base $M$ such that 
$SC_x (M) \subset T_x (M), \;\; x \in M$. The pair 
$(M, SC (p, q))$ is said to be  an {\it almost Grassmann 
structure} and is denoted by $AG (p-1, p+q-1)$. The manifold 
$M$ endowed with  such a structure is said to be an {\it almost 
Grassmann manifold}.
}

As was the case for Grassmann structures, the almost Grassmann 
structure $AG (p-1, p+q-1)$ is equivalent to the structure 
$AG (q-1, p+q-1)$ since both of these structures are generated on 
the manifold $M$ by a differentiable family of  Segre cones 
$SC_x (p, q)$.

 Let us consider some examples. 

\examp{\label{examp:2.2}
The main example of an almost Grassmann structure is the almost 
Grassmann structure associated with the Grassmannian $G (m, n)$.  
As we saw, there is a field of Segre cones $SC_x (p, q) 
= T_x (\Omega) \bigcap \Omega, \; x \in \Omega$, where $p = m+1$ 
and $q = n-m$, which defines  an almost Grassmann structure.
}
 
\examp{\label{examp:2.3}
Consider a pseudoconformal $CO (2, 2)$-structure on a 
four-di\-men\-sional manifold $M$. The isotropic cones $C_x$ of 
this structure carry two families of plane generators. Hence, 
these cones are Segre cones $SC_x (2, 2)$. Therefore, a 
pseudoconformal $CO (2, 2)$-structure is an almost Grassmann 
structure $AG (1, 3)$. 

If we complexify the four-dimensional tangent subspace 
$T_x (M^4)$ and consider Segre cones with complex generators, 
then  conformal $CO (1, 3)$- and $CO (4, 0)$-structures can also 
be considered as complex almost Grassmann structures of the same 
type $AG (1, 3)$. However, in this paper, we will consider only 
real almost Grassmann structures.
}

Almost Grassmann structures arise also in the study of 
multidimensional webs (see  \cite{AS92} and \cite{G88}). 

\examp{\label{examp:2.4}
Consider a three-web formed on a manifold $M^{2q}$ of dimension 
$2q$ by three foliations $\lambda_u, \; u = 1, 2, 3$, of 
codimension $q$ which are in general position (see \cite{A69} 
or \cite{AS92}).   
Through any point $x \in M^{2q}$, there pass three leaves 
${\cal F}_u$ belonging to the foliations $\lambda_u$. In the 
tangent subspace $T_x (M^{2q})$, we consider three subspaces 
$T_x ({\cal F}_u)$ which are tangent to ${\cal F}_u$ at the point 
$x$. If we take the projectivization of this configuration with  
center at the point $x$, then we obtain a projective space 
$P^{2q-1}$ of dimension $2q-1$ containing three subspaces of 
dimension $q-1$ which are in general position. These three 
subspaces determine a Segre variety $S (1, q - 1)$, and the 
latter variety is the directrix for a Segre cone 
$SC_x (2, q) \subset T_x (M^{2q})$. Thus, on $M^{2q}$, a field of 
Segre cones  arises, and this field determines an almost 
Grassmann structure on $M^{2q}$. 

The structural group   of the web $W (3, 2, q)$ 
is smaller than that  of the induced  almost Grassmann structure, 
since transformations of this group must keep invariant the 
subspaces $T_x ({\cal F}_u)$. Thus, the structural group of 
the three-web is the group ${\bf GL} (q)$. 
}

\examp{\label{examp:2.5} 
 Consider a $(p+1)$-web $W(p+1, p, q) 
= (M; \lambda_1, \ldots, \lambda_{p+1})$ formed on a 
differentiable manifold $M$ of dimension $pq$ by $p+1$ foliations 
$\lambda_u, \; u=1, \ldots, p+1$ of dimension $q$ which are in 
general position on $M$ (see \cite{G73} or \cite{G88}).

As in Example {\bf 2.4}, the tangent spaces $T_x ({\cal F}_u)$ 
define the cone $SC_x (p, q) \supset T_x ({\cal F}_u)$, and 
the field of these cones defines an almost Grassmann structure 
$AG (p-1, p+q-1)$ on $M$.

The structural group of the web $W (p+1, p, q)$ is the same group 
 $G ={\bf GL} (q)$ as for the web $W (3, 2, q)$, and 
this group does not depend on $p$.
}

{\bf 2.} 
The structural group of the almost Grassmann structure is a 
subgroup of the general linear group ${\bf GL}(pq)$ of 
transformations of the space $T_{x} (M)$, which leave the cone 
$SC_{x}(p, q) \subset T_x (M)$ invariant. We denote this group by 
$G = {\bf GL} (p, q)$. 

To clarify the structure of this group, in the tangent space  
$T_x (M)$, we consider a family of frames  
$\{e_{i}^{\alpha}\}, \alpha = 1, \ldots, p; \;\; i=p+1, \ldots, 
p+q$, such that for any fixed $i$, the vectors $e_i^\alpha$ 
belong to a $p$-dimensional generator $\xi$ of the Segre cone 
$SC_x (p, q)$, and for any fixed $\alpha$, the vectors 
$e_i^\alpha$ belong to a $q$-dimensional generator 
$\eta$ of  $SC_x (p, q)$.  In such a frame, the equations of 
the cone $SC_x (p, q)$ can be written as follows:
\begin{equation}\label{eq:2.1}
z_\alpha^i = t_\alpha s^i, \;\; \alpha = 1, \ldots, p, \;\; 
i = p +1, \ldots, p + q,
\end{equation}
where $z^i_\alpha$ are the coordinates of a vector 
$z = z^i_\alpha e_i^\alpha \subset T_x (M)$, and $t_\alpha$ and 
$s^i$ are parameters on which a vector $z \subset SC_x (M)$ 
depends. 

The family of frames $\{e_i^\alpha\}$ attached to the cone 
$SC_x (p, q)$ admits a transformation of the form 
\begin{equation}\label{eq:2.2}
'e_i^\alpha = A_\beta^\alpha A_i^j e_j^\beta,
\end{equation}
where $(A_\beta^\alpha)$ and $(A_j^i)$ are nonsingular square 
matrices of orders $p$ and $q$, respectively. These matrices are 
not defined uniquely since they admit a multiplication by  
reciprocal scalars. However, they can be made unique by 
restricting to unimodular matrices  $(A_\beta^\alpha)$ or 
$(A_i^j)$: $\det (A_\beta^\alpha) = 1$ or  $\det (A_i^j) = 1$. 
Thus the structural group of the almost Grassmann structure 
defined by equations (2.2), can be represented in the form 
\begin{equation}\label{eq:2.3}
G = {\bf SL} (p) \times {\bf GL}(q) \cong {\bf GL} (p) \times 
{\bf SL} (q),
\end{equation}
where ${\bf SL} (p)$ and ${\bf SL} (q)$ are special linear groups 
of dimension $p$ and $q$, respectively. Such a representation has 
been used by  Hangan (\cite{H66}, \cite{H68}, \cite{H80}) 
and  Goldberg \cite{G75} (see also 
the book \cite{G88}, Ch. 2), and   \cite{M78}. Unlike this 
approach, we will assume that both matrices  
$(A_{\beta}^{\alpha})$ and  $(A_i^j)$ are unimodular but 
the right-hand side of equation (2.2) admits a multiplication 
by a scalar factor. As a result, we obtain a more symmetric 
representation of the group $G$:
\begin{equation}\label{eq:2.4}
G = {\bf SL} (p) \times {\bf SL}(q)  \times {\bf H},
\end{equation}
where ${\bf H} = {\bf R}^* \otimes \mbox{{\rm Id}}$ is the group  
of homotheties of the $T_x (M)$. 

It follows that {\em an almost Grassmann structure  $AG (m, n)$ 
is a $G$-structure of first order.}

It follows from condition (2.1) that $p$-dimensional plane 
generators $\xi$ of the Segre cone $SC_x (p, q)$ are determined 
by values of the parameters $s^i$, and $t_\alpha$ are coordinates 
of points of a generator $\xi$. But a plane generator $\xi$ is 
not changed if we multiply the parameters $s^i$ by the same 
number. Thus, the family of plane generators $\xi$ depends on 
$q - 1$ parameters.

Similarly, $q$-dimensional plane generators $\eta$ of the Segre 
cone $SC_x (p, q)$ are determined by values of the parameters 
$t_\alpha$, and $s^i$ are coordinates of points of a generator 
$\eta$. 
But a plane generator $\eta$ is not changed if we multiply 
the parameters $t_\alpha$ by the same number. Thus, the family 
of plane generators $\eta$ depends on $p - 1$ parameters.

The $p$-dimensional subspaces $\xi$ form a fiber bundle on 
the manifold $M$. The base of this bundle is the manifold $M$, 
and its fiber attached to a point $x \in M$ is the set of all 
$p$-dimensional plane generators $\xi$ of the Segre cone 
$SC_x (p, q)$. The dimension of a fiber is $q - 1$, and it is 
parametrized by means of a projective space 
$P_\alpha, \; \dim P_\alpha = q - 1$. We will denote this fiber 
bundle of $p$-subspaces by $E_\alpha = (M, P_\alpha)$. 

In a similar manner, $q$-dimensional plane generators $\eta$ of 
the Segre cone $SC_x (p, q)$ form on $M$ the fiber bundle 
$E_\beta = (M, P_\beta)$ with the base $M$ and fibers of 
dimension $p - 1 = \dim P_\beta$. The fibers are $q$-dimensional 
plane generators $\eta$ of the Segre cone $SC_x (p, q)$. 
 
Consider the manifold $M_\alpha = M \times P_\alpha$ of dimension 
$pq + q - 1$. The fiber bundle $E_\alpha$ induces on $M_\alpha$ 
the distribution $\Delta_\alpha$ of plane elements $\xi_\alpha$ 
of dimension $q$. In a similar manner, on the manifold 
$M_\beta = M \times P_\beta$ the fiber bundle $E_\beta$ induces 
the distribution $\Delta_\beta$ of plane elements $\eta_\beta$ of 
dimension $p$. 

\defin{\label{def:2.6}
An almost Grassmann structure $AG (p-1, p+q-1)$ is said to be 
{\em $\alpha$-semiintegrable} if the distribution $\Delta_\alpha$ 
is  integrable on this structure.  Similarly, an almost Grassmann 
structure $AG (p-1, p+q-1)$ is said to be 
{\em $\beta$-semiintegrable} if  the distribution $\Delta_\beta$ 
is  integrable on this structure. A structure $AG (p-1, p+q-1)$ 
is called {\em integrable} if it is both $\alpha$- and 
$\beta$-semiintegrable.
} 

Integral manifolds $\widetilde{V}_\alpha$ of the distribution 
$\Delta_\alpha$ of an $\alpha$-semiintegrable almost Grassmann 
structure are of dimension $p$. They are projected on the 
original manifold $M$ in the form of a submanifold $V_\alpha$ of 
the same dimension $p$, which, at any of its points, is tangent 
to the $p$-subspace $\xi_\alpha $ of the fiber bundle $E_\alpha$. 
Through each point $x \in M$, there passes a $(q-1)$-parameter 
family of submanifolds $V_\alpha$.

Similarly, integral manifolds $\widetilde{V}_\beta$ of the 
distribution $\Delta_\beta$ of a $\beta$-semiinte\-gra\-ble 
almost 
Grassmann structure are of dimension $q$. They are projected on 
the original manifold $M$ in the form of a submanifold $V_\beta$ 
of the same dimension $q$, which, at any of its points, is 
tangent to the $q$-subspace $\eta_\beta$ of the fiber bundle 
$E_\beta$. Through each point $x \in M$, there passes a 
$(p-1)$-parameter family of submanifolds $V_\beta$. 

If an almost Grassmann structure on $M$ is integrable, then 
through each point $x \in M$, there pass a $(q-1)$-parameter 
family of submanifolds $V_\alpha$ and a $(p-1)$-parameter family 
of submanifolds $V_\beta$ which were described above.

The Grassmann structure $G (m, n)$ is an integrable almost 
Grassmann structure $AG (m, n)$ since through 
any point $x \in \Omega (m, n)$, onto which the manifold 
$G (m, n)$ is mapped bijectively  under the Grassmann mapping, 
there pass a $(q - 1)$-parameter family of $p$-dimensional plane 
generators (which are the submanifolds $V_\alpha$) and a 
$(p - 1)$-parameter family of $q$-di\-men\-sional 
plane generators 
(which are the submanifolds $V_\beta$). In the projective space 
$P^n$, there corresponds to  submanifolds $V_\alpha$ a family of 
$m$-dimensional subspaces belonging to a subspace of dimension 
$m + 1$, and  there corresponds to  submanifolds $V_\beta$ a 
family of $m$-dimensional subspaces passing through a subspace of 
dimension $m - 1$.

\section{Structure Equations of an Almost Grassmann \newline 
Structures}

\setcounter{equation}{0}
 
{\bf 1.} 
Consider  a differentiable manifold $M$ of dimension $pq$ 
endowed with an almost Grassmann structure 
$AG (p-1, p+q-1)$. Suppose that $x \in M$,  $T_x (M)$ 
is the tangent space of the manifold    $M$ at the 
point $x$ and that $\{e_i^\alpha\}$ is an adapted frame 
of the structure $AG (p-1, p+q-1)$. The decomposition of 
a vector $z \in T_x (M)$ with respect to this basis can be 
written in the form
$$
z = \omega_\alpha^i (z) e_i^\alpha,
$$
where $\omega_\alpha^i$ are 1-forms making up the {\em co-frame} 
in the space $T_x (M)$. If $z = dx$ is the differential of a point $x \in M$, then the forms $\omega_\alpha^i (dx)$ are 
differential forms defined on a first-order frame bundle 
associated with the almost Grassmann structure. These forms 
constitute a completely integrable system of forms. As a 
result we have
\begin{equation}\label{eq:3.1}
 d\omega_\alpha^i = \omega_\beta^j \wedge 
\omega_{\alpha j}^{i\beta}.
\end{equation}
The forms $\omega_\alpha^i$ are called also the {\em basis forms} 
of the manifold $M$.

As earlier, we set $\pi_{\beta i}^{j \alpha} 
= \omega_{\beta i}^{j \alpha}  (\delta)$, where $\delta$ 
is the operator of differentiation with respect to 
the fiber parameters of the frame bundle. These forms 
determine an infinitesimal transformation of the adapted frames:
\begin{equation}\label{eq:3.2}
 \delta e_i^\alpha = \pi_{\beta i}^{j \alpha}   e_j^\beta.
\end{equation}
On the other hand, the admissible transformations of 
adapted frames can be written as closed form equations (2.2). 
Solving equations (2.2), we obtain
\begin{equation}\label{eq:3.3}
e_{i}^{\alpha}= \tilde{A}_{\beta}^{\alpha} \tilde{A}_{i}^{j} \:
 'e_{j}^{\beta},
\end{equation}
where $(\tilde{A}_{\beta}^{\alpha})$ and $(\tilde{A}_{j}^{i})$ are the inverse matrices 
 of the matrices $(A_{\beta}^{\alpha})$ and $(A_{j}^{i})$, 
respectively: 
\begin{equation}\label{eq:3.4}
A_{\gamma}^{\alpha}  \tilde{A}_{\beta}^{\gamma} = 
A_{\beta}^{\gamma}\tilde{A}_{\gamma}^{\alpha}
= \delta_{\beta}^{\alpha},\;\;\;
A_k^i \tilde{A}_{j}^{k}=A_{j}^{k} \tilde{A}_{k}^{i}=\delta_{j}^{i}.
\end{equation}
It follows from (3.4) that
\begin{equation}\label{eq:3.5}
A^{\gamma}_{\beta} \cdot \delta \tilde{A}^{\alpha}_{\gamma} = 
-\tilde{A}_{\gamma}^{\alpha} \cdot \delta A^{\gamma}_{\beta},\;\;
A^{i}_{k} \cdot \delta \tilde{A}_{j}^{k} 
= - \tilde{A}_{j}^{k} \cdot \delta A_{k}^{i},
\end{equation}

Suppose now that $\{x,\;'e_i^\alpha\}$ is a fixed frame, 
 $d('e_i^\alpha) = 0$. Then, differentiating (3.3) for 
a fixed $x \in M$ and using (3.4) and (3.5), we obtain
\begin{equation}\label{eq:3.6}
 \delta e_i^\alpha=(\delta_\beta^\alpha \pi_i^j -
 \delta_i^j \pi_\beta^\alpha) e_j^\beta,
\end{equation}
where
\begin{equation}\label{eq:3.7}
 \pi_i^j=A_k^j \cdot \delta \tilde{A}_i^k,\;\;
\pi_\beta^\alpha= \tilde{A}_\gamma^\alpha \cdot \delta 
A_\beta^\gamma.
\end{equation}
Comparing formulas  (3.2) and (3.6), we find that   
\begin{equation}\label{eq:3.8}
 \pi_{\beta i}^{j\alpha} 
=\delta^{\alpha}_{\beta}\pi_{i}^{j} 
- \delta_{i}^{j}\pi_{\beta}^{\alpha}.
\end{equation}
In these formulas the forms $\pi_j^i$ are 
invariant forms of the group ${\bf GL} (q)$, 
the forms $\pi_\beta^\alpha$ are 
invariant forms of the group ${\bf GL} (p)$, and 
the forms $\pi_{\beta i}^{j\alpha}$ are 
invariant forms of the structural 
group  $G$ of the almost Grassmann 
structure  $AG (p-1, p+q-1)$. 

If a point $x \in M$ is  variable, then from equations 
(3.8) we find that
\begin{equation}\label{eq:3.9}
 \omega_{\beta i}^{j\alpha} 
= \delta^{\alpha}_{\beta}\omega_{i}^{j} - 
 \delta_{i}^{j}\omega_{\beta}^{\alpha}+ \widehat{u}_{\beta ik}^{j\alpha\gamma} \omega_{\gamma}^{k},
\end{equation}
where $\widehat{u}_{\beta ik}^{j\alpha\gamma}$ are certain 
functions  defined on the first-order frame bundle. 
Substituting for $\omega_{\beta i}^{j\alpha}$ in (3.1)
 their values taken   from (3.9), we obtain
\begin{equation}\label{eq:3.10}
 d\omega_\alpha^i = 
 \omega_\alpha^\beta \wedge \omega_\beta^i 
+  \omega_\alpha^j \wedge \omega_j^i 
+ u_{\alpha jk}^{i\beta\gamma} 
\omega_\beta^j \wedge \omega_\gamma^k,
\end{equation}
where $u_{\alpha jk}^{i\beta\gamma}$ denotes the result of 
alternation of the quantities 
$\widehat{u}_{\alpha jk}^{i\beta\gamma}$ occurring in (3.9) 
 with respect to the pairs of indices ${\beta \choose j}$ and 
 ${\gamma \choose k}$: 
$ u_{\alpha jk}^{i\beta\gamma}= -u_{\alpha kj}^{i\gamma\beta}.$

If we set 
$$
\omega_\alpha^\beta = \widetilde{\omega}_\alpha^\beta 
+ \frac{1}{p} \delta_\alpha^\beta \omega^\gamma_\gamma, \;\;
\omega_i^j = \widetilde{\omega}_i^j  
+ \frac{1}{q} \delta_i^j  \omega^k_k, 
$$
then it is easy to see that $\widetilde{\omega}_\alpha^\alpha 
= 0$ and $\widetilde{\omega}_k^k = 0$, and the above structure 
equations take the form
$$
 d\omega_{\alpha}^{i}=
 \widetilde{\omega}_{\alpha}^{\beta} \wedge \omega_{\beta}^{i}  
 + \omega_{\alpha}^{j} \wedge \widetilde{\omega}_{j}^{i} 
+ \omega \wedge \omega_\alpha^i + 
u_{\alpha jk}^{i\beta\gamma} 
\omega_{\beta}^{j} \wedge \omega_{\gamma}^{k},
$$
where 
$\omega = 
\frac{1}{p} \omega^\gamma_\gamma - 
\frac{1}{q}   \omega^k_k$. If we suppress 
$\sim$, then the structure equations take the form:
\begin{equation}\label{eq:3.11}
 d\omega_\alpha^i = 
\omega_\alpha^\beta \wedge \omega_\beta^i 
 + \omega_\alpha^j \wedge \omega_j^i 
+  \omega \wedge \omega_\alpha^i + u_{\alpha jk}^{i\beta\gamma} 
\omega_{\beta}^{j} \wedge \omega_{\gamma}^{k},
\end{equation}
where 
\begin{equation}\label{eq:3.12}
 u_{\alpha jk}^{i\beta\gamma}= - u_{\alpha kj}^{i\gamma\beta}
\end{equation}
and 
\begin{equation}\label{eq:3.13}
\omega_\gamma^\gamma = 0, \;\; \omega_k^k = 0. 
\end{equation}

Conditions (3.13) mean that the subgroups 
${\bf GL} (p)$ and ${\bf GL} (q)$ 
of the structural group $G$ of the almost Grassmann structure 
$AG (p - 1, p + q - 1)$ are reduced to the groups 
${\bf SL} (p)$ and ${\bf SL} (p)$, respectively, 
and that the group $G$ itself is represented in the form 
(2.4). As for the Grassmannian $G (p - 1, p + q - 1)$ 
(see Section {\bf 1}), for the almost Grassmann manifold 
the forms $\omega_\alpha^\beta, \omega_i^j$, and $\omega$ are 
fiber forms defined on the second-order frame bundle associated 
with the almost Grassmann manifold $AG (p-1, p+q-1)$. 

The structure equations (3.11) differ from the structure 
equations (1.7) of the Grassmannian  $G (p - 1, p + q - 1)$ 
only by the last term.

{\bf 2.} 
We obtain the remaining structure equations of the almost 
Grassmann manifold $M$ by exterior differentiation of
 (3.11). This gives 

\begin{equation}\label{eq:3.14}
\renewcommand{\arraystretch}{1.5}
\begin{array}{ll}
  \Omega_{\alpha}^{\beta} \wedge \omega_{\beta}^{i} 
 - \Omega_{j}^{i} \wedge \omega_{\alpha}^{j}
+ (\nabla u_{\alpha jk}^{i\beta\gamma} 
+  u_{\alpha jk}^{i\beta\gamma} \omega)
 \wedge \omega_{\beta}^{j} \wedge \omega_{\gamma}^{k} \\
+ d \omega \wedge \omega_\alpha^i 
+ 2u_{\alpha mk}^{i\varepsilon\gamma} 
u_{\varepsilon lj}^{m\delta\beta} 
 \omega_{\delta}^{l} \wedge \omega_{\beta}^{j} \wedge 
\omega_{\gamma}^{k} = 0,
\end{array}
\renewcommand{\arraystretch}{1}
\end{equation}
where
$$
\renewcommand{\arraystretch}{1.5}
\begin{array}{ll}
\Omega_{\alpha}^{\beta}  =  d\omega_{\alpha}^{\beta} 
- \omega_{\alpha}^{\gamma}\wedge \omega_{\gamma}^{\beta}, \;\;
\Omega_{j}^{i}  =  d\omega_{j}^{i} - \omega_{j}^{k}
\wedge \omega_{k}^{i}, \\ 
\nabla u_{\alpha jk}^{i\beta\gamma} 
= du_{\alpha jk}^{i\beta\gamma}
- u_{\delta jk}^{i\beta\gamma} \omega_{\alpha}^{\delta} 
- u_{\alpha lk}^{i\beta\gamma}\omega_{j}^{l} 
- u_{\alpha jl}^{i\beta\gamma}\omega_{k}^{l} 
  + u_{\alpha jk}^{l\beta\gamma}\omega_{l}^{i}
+ u_{\alpha jk}^{i\delta\gamma}\omega_{\delta}^{\beta}
+ u_{\alpha jk}^{i\beta\delta}\omega_{\delta}^{\gamma}.
\end{array}
\renewcommand{\arraystretch}{1}
$$

To solve equations (3.14), we represent the 
forms $\Omega_{j}^{i},  \Omega_{\alpha}^{\beta}$, 
and $d \omega$ 
 from the left-hand side of equations 
(3.14) as a sum of terms containing the basis forms 
and the terms not containing these forms:
\begin{equation}\label{eq:3.15}
 \Omega_\alpha^\beta
=  \omega_{\gamma}^{k} \wedge \omega_{\alpha k}^{\beta\gamma} 
+  \Phi_{\alpha}^{\beta}, \;\;
 \Omega_j^i 
=  \omega_{jk}^{i \gamma} \wedge \omega_{\gamma}^{k} 
+ \Phi_{j}^{i},\;\;
d \omega = \omega_i^\alpha \wedge \omega_\alpha^i + \Phi,
\end{equation}
where $\Phi_\alpha^\beta, \Phi_j^i, \Phi$ and  
$\omega_{\alpha k}^{\beta \gamma}, \omega_{jk}^{i\gamma}$ are 
certain 2- and 1-forms not expressed in terms of 
 the basis forms $\omega_{\alpha}^{i}$ only. 

By (3.13), we have $\Omega^\gamma_\gamma = 0$ 
and $\Omega_k^k = 0$, which implies that 
\begin{equation}\label{eq:3.16}
  \Phi_\gamma^\gamma =0, \;\; \Phi_k^k = 0.
\end{equation}
and 
\begin{equation}\label{eq:3.17}
\omega_{\alpha k}^{\alpha \gamma} = 0, \;\; 
\omega_{ik}^{i\gamma} = 0 
\end{equation}

 Substituting (3.15)  into equations (3.14), 
 we obtain
\begin{equation}\label{eq:3.18}
\renewcommand{\arraystretch}{1.5}
 \begin{array}{ll}
(\delta_j^i  \Phi_\alpha^\beta - \delta_\alpha^\beta \Phi_j^i 
+ \delta_\alpha^\beta \delta_j^i \Phi) \wedge \omega_\beta^j 
+ 2u_{\alpha m[k}^{i\varepsilon[\gamma} 
u_{|\varepsilon| lj]}^{|m|\delta\beta]}
 \omega_\delta^l \wedge \omega_\beta^j \wedge 
\omega_\gamma^k \\
+ (\delta_{[j}^i  \omega_{|\alpha| k]}^{[\beta \gamma]} 
+ \delta_\alpha^{[\beta} \omega_{[jk]}^{|i|\gamma]}  
+ \delta_\alpha^{[\gamma} \delta_{[k}^{|i|} \omega_{j]}^{\beta]} 
+ \nabla u^{i\beta\gamma}_{\alpha jk} 
+  u^{i\beta\gamma}_{\alpha jk} \omega) \wedge 
\omega_\beta^j \wedge \omega_\gamma^k = 0,
 \end{array} 
\renewcommand{\arraystretch}{1}
\end{equation}
where the alternation is carried over with respect to vertical pairs 
of indices.
The first term in the left-hand side of (3.18) 
does not have similar terms among other terms of this side. 
Thus this term vanishes. But since the first factor of this term 
does not contain the basis forms, this factor itself vanishes: 
\begin{equation}\label{eq:3.19}
\delta_j^i  \Phi_\alpha^\beta - \delta_\alpha^\beta \Phi_j^i 
+ \delta_\alpha^\beta \delta_j^i \Phi = 0.
\end{equation}
Contracting (3.19) with respect to the indices $\alpha$ and 
$\beta$, applying (3.16),  and dividing by $p$, we find that 
$$
- \Phi_j^i + \delta_j^i \Phi = 0.
$$
Contracting this equation with respect to the indices $i$ 
and $j$, we obtain $\Phi = 0$, and consequently $\Phi^i_j = 0$. 
Finally, by (3.19), we find that $\Phi_\alpha^\beta = 0$.

Now equation (3.18) contains only the last 
two terms. It follows that the 1-form which is multiplied 
by $\omega_\beta^j \wedge \omega_\gamma^k $ is expressed only 
in terms of the basis forms. Therefore, if the principal 
parameters are fixed (i.e., if $\omega_\alpha^i = 0$), 
then we obtain
\begin{equation}\label{eq:3.20}
 2 (\nabla_\delta u^{i\beta\gamma}_{\alpha jk} 
+  u^{i\beta\gamma}_{\alpha jk} \pi) 
+ \delta_j^i  \pi_{\alpha k}^{\beta \gamma} 
- \delta_k^i  \pi_{\alpha j}^{\gamma \beta} 
+ \delta_\alpha^\beta \pi_{jk}^{i\gamma}  
- \delta_\alpha^\gamma \pi_{kj}^{i\beta}  
+ \delta_\alpha^{\gamma} \delta_k^i \pi_j^\beta 
- \delta_\alpha^{\beta} \delta_j^i \pi_k^\gamma 
= 0, 
\end{equation}
where as usual $\pi_k^\gamma = \omega_k^\gamma (\delta), \;
\pi_{jk}^{\beta\gamma} = \omega_{jk}^{\beta\gamma} (\delta)$, 
and $\pi = \omega (\delta)$.
It follows from equation (3.20) that the quantities 
$u^{i\beta\gamma}_{\alpha jk}$ form a geometric object 
that is defined in a second-order differential neighborhood of 
the almost Grassmann structure $AG (p - 1, p + q - 1)$.

Consider the quantities 
\begin{equation}\label{eq:3.21}
u_{\alpha k}^{\beta \gamma} = u_{\alpha ik}^{i\beta \gamma}, 
\;\;u_{j k}^{i \gamma} = u_{\alpha jk}^{i \alpha\gamma}.
\end{equation}
If we contract  equations  (3.20) with respect to 
the indices $i$ and $j$, then after some calculations we find 
that 
\begin{equation}\label{eq:3.22}
\nabla_\delta u_{\alpha k}^{\beta \gamma} + u_{\alpha k}^{\beta 
\gamma} \pi = - \displaystyle \frac{1}{2} \Bigl[q (\pi_{\alpha 
k}^{\beta \gamma} - \delta_{\alpha}^{\beta} \pi_{k}^{\gamma}) 
- \pi_{\alpha k}^{\gamma \beta} - \delta_{\alpha}^{\gamma} 
(\pi_{ki}^{i \beta} - \pi_{k}^{\beta})\Bigr].
\end{equation}
 Similarly, 
contracting equations (3.20) with respect to 
the indices $\alpha$ and $\beta$, we obtain
\begin{equation}\label{eq:3.23}
\nabla_\delta u_{j k}^{i \gamma} + u_{j k}^{i \gamma} \pi = 
- \displaystyle \frac{1}{2} \Bigl[p (\pi_{jk}^{i \gamma} 
- \delta_j^i \pi_k^\gamma) - \pi_{kj}^{i \gamma} - 
\delta_k^i (\pi_{\alpha j}^{\gamma \alpha} 
- \pi_{j}^{\gamma})\Bigr].
\end{equation}
Formulas (3.22) and (3.23) show that each of the quantities 
$u_{\alpha k}^{\beta \gamma}$ and  $u_{j k}^{i \gamma}$ form a geometric object that is defined in a second-order 
differential neighborhood of 
the almost Grassmann structure $AG (p - 1, p + q - 1)$. 

Let us prove that if we make 
a specialization of second-order frames, then we can reduce 
these geometric objects to 0. 

We will prove this for the geometric object 
$u_{\alpha k}^{\beta \gamma}$. To this end, we must show that the 
1-forms in the right-hand sides of equations (3.22) are 
linearly independent. First, we note that the forms 
$\pi_{\alpha k}^{\beta \gamma}$ are linearly 
independent in the set of  second-order frames. Let us equate 
to 0 the right-hand sides of equations (3.22):
\begin{equation}\label{eq:3.24}
q (\pi_{\alpha k}^{\beta \gamma} 
- \delta_{\alpha}^{\beta} \pi_{k}^{\gamma}) 
- \pi_{\alpha k}^{\gamma \beta} - \delta_{\alpha}^{\gamma} 
(\pi_{ki}^{i \beta} - \pi_{k}^{\beta}) =0.
\end{equation}
If we contract equations (3.24) first 
 with respect to the indices $\alpha$ and $\beta$ and second 
 with respect to the indices $\alpha$ and $\gamma$, we arrive at 
the system
\begin{equation}\label{eq:3.25}
\renewcommand{\arraystretch}{1.5}
\left\{
 \begin{array}{ll}
\pi_{\alpha k}^{\gamma \alpha} +  \pi_{ki}^{i\gamma} 
= (1 - pq) \pi_{k}^{\gamma}, \\
q \pi_{\alpha k}^{\gamma\alpha} - p \pi_{ki}^{i\gamma} 
= (q - p) \pi_k^\gamma.
 \end{array} 
\right.
\renewcommand{\arraystretch}{1}
\end{equation}
If we solve this system, we find the quantities $\pi_{\alpha k}^{\gamma \alpha}$ and $\pi_{ki}^{i\gamma}$:
\begin{equation}\label{eq:3.26}
\pi_{\alpha k}^{\gamma \alpha} = - \frac{q (p^2 - 1)}{p+q} 
 \pi_{k}^{\gamma}, \;\;
\pi_{ki}^{i \gamma} = - \frac{p (q^2 - 1)}{p+q} 
 \pi_{k}^{\gamma}.
\end{equation}
Substituting these values of $\pi_{ki}^{i\gamma}$ 
into equations (3.24), after some calculations, we 
reduce the equations obtained to the following form:
\begin{equation}\label{eq:3.27}
q \widetilde{\pi}_{\alpha k}^{\beta \gamma} 
- \widetilde{\pi}_{\alpha k}^{\gamma \beta} =0,
\end{equation}
where
$$
\widetilde{\pi}_{\alpha k}^{\beta \gamma} = 
\pi_{\alpha k}^{\beta \gamma} - \displaystyle \frac{q}{p+q} 
\Bigl(\delta_{\alpha}^{\beta} \pi_{k}^{\gamma} - 
p \delta_{\alpha}^{\gamma} \pi_{k}^{\beta}\Bigr).
$$
Interchanging in (3.27) the indices $\beta$ and $\gamma$, 
we obtain
\begin{equation}\label{eq:3.28}
- \widetilde{\pi}_{\alpha k}^{\beta \gamma} 
+ q \widetilde{\pi}_{\alpha k}^{\gamma \beta} = 0.
\end{equation}
Since the determinant of the system of equations 
(3.27)--(3.28) is equal to $q^2 - 1 \neq 0$, the 
system has only the trivial solution.

But the forms $\widetilde{\pi}_{\alpha k}^{\beta \gamma} $ as 
the forms $\pi_{\alpha k}^{\beta \gamma}$ are linearly 
independent. Thus the forms 
$q \widetilde{\pi}_{\alpha k}^{\beta \gamma} 
- \widetilde{\pi}_{\alpha k}^{\gamma \beta}$ are linearly 
independent too. But, up to the factor 
$- 
\frac{1}{2}$, the latter forms coincide with 
the right-hand sides of equations (3.22). 

Hence the geometric object 
$u_{\alpha k}^{\beta \gamma} = u_{\alpha ik}^{i\beta \gamma}$ 
can be reduced to 0. Similarly the geometric object 
$u_{j k}^{i \gamma} = u_{\alpha jk}^{i\alpha \gamma}$ 
can be reduced to 0. This operation leads to a reduction 
of the set of second-order frames of 
the almost Grassmann structure $AG (p - 1, p + q - 1)$. 
Before this reduction, the set of second-order frames 
depended on $pq (p^2 + q^2)$ parameters  equal 
to the number of linearly independent forms among the forms 
$\pi_{\alpha k}^{\beta \gamma}$ and $\pi_{jk}^{i \gamma}$. 
After the reduction, the forms 
$\widetilde{\pi}_{\alpha k}^{\beta \gamma}$ and 
$\widetilde{\pi}_{jk}^{i \gamma}$ vanish, 
and the forms $\pi_{\alpha k}^{\beta \gamma}$ and 
$\pi_{jk}^{i \gamma}$ are expressed 
in terms of the 1-forms $\pi_{k}^{\gamma}$: 
\begin{equation}\label{eq:3.29}
\pi_{\alpha k}^{\beta \gamma} 
= \displaystyle \frac{q}{p+q} \Bigl(\delta_\alpha^\beta 
 \pi_{k}^{\gamma} - p \delta_\alpha^\gamma \pi_{k}^{\beta}\Bigr), 
\;\; 
\pi_{jk}^{i \gamma} 
= \displaystyle \frac{p}{p+q} \Bigl(\delta_j^i 
 \pi_k^\gamma - q \delta_k^i \pi_j^\gamma\Bigr). 
\end{equation}
Since there are $pq$ forms $\pi_k^\gamma$, and they are linearly 
independent, the reduced family of second-order frames depends on 
$pq$ parameters. The 1-forms $\pi^\gamma_k$ define admissible 
transformations of frames in this reduced family of second-order 
frames.

Denote by $a^{i\beta\gamma}_{\alpha jk}$ 
 the quantities $u^{i\beta\gamma}_{\alpha jk}$ after the 
specialization indicated above. 
Then the quantities $a^{i\beta\gamma}_{\alpha jk}$ satisfy 
the conditions
\begin{equation}\label{eq:3.30}
a^{i\alpha\gamma}_{\alpha jk} = 0, \;\;\;
a^{i\beta\gamma}_{\alpha ik} = 0
\end{equation}
and
\begin{equation}\label{eq:3.31}
a^{i\beta\gamma}_{\alpha jk} = - a^{i\gamma\beta}_{\alpha kj}.
\end{equation}
The last relations follow from conditions (3.12). 

Substituting expressions (3.29) into equations (3.20), 
we find that 
\begin{equation}\label{eq:3.32}
\nabla_\delta a^{i\beta\gamma}_{\alpha jk} + 
a^{i\beta\gamma}_{\alpha jk} \pi = 0.
\end{equation}
This implies the following theorem: 

\begin{theorem}
 The quantities $a^{i\beta\gamma}_{\alpha jk}$, defined in a 
second-order neighborhood by the reduction of second-order frames 
indicated above, form a relative tensor of weight $- 1$ 
and satisfy conditions $(3.30)$ and $(3.31)$.
\end{theorem}

\defin{\label{def:3.3}
The tensor $\{a_{\alpha jk}^{i\beta\gamma}\}$ is said to be the 
{\it  first structure tensor}, or the {\it torsion tensor}, of an 
almost Grassmann manifold $AG (p-1, p+q-1)$.
}

After the specialization of second-order frames 
has been made, the first structure equations (3.11) 
become
\begin{equation}\label{eq:3.33}
 d\omega_\alpha^i = \omega_\alpha^j \wedge 
\omega_j^i + \omega_\alpha^\beta \wedge \omega_\beta^i 
+  \omega \wedge \omega_\alpha^i + a_{\alpha jk}^{i\beta\gamma} 
\omega_{\beta}^{j} \wedge \omega_{\gamma}^{k}.
\end{equation}

{\bf 3}. We will now find the expression for the tensor 
$a_{\alpha jk}^{i\beta\gamma}$ in terms of the quantities 
$u_{\alpha jk}^{i\beta\gamma}$ occurring in equations 
(3.11). We assume that the specialization of 
second-order frames 
indicated above has not been made and that the quantities 
$u_{\alpha jk}^{i\beta\gamma}$ satisfy equations (3.20), 
which we write in the form 
\begin{equation}\label{eq:3.34}
\renewcommand{\arraystretch}{1.5}
 \begin{array}{ll}
\nabla_\delta u^{i\beta\gamma}_{\alpha jk} 
+  u^{i\beta\gamma}_{\alpha jk} \pi 
=& \!\!\!\!\displaystyle \frac{1}{2} \Bigl(\delta_\alpha^{\gamma} \pi_{kj}^{i\beta} 
- \delta_j^i  \pi_{\alpha k}^{\beta \gamma} 
+ \delta_k^i \pi_{\alpha j}^{\gamma \beta} 
- \delta_\alpha^\beta \pi_{jk}^{i\gamma}\\
&\!\!\!\!+ \delta_j^i \delta_\alpha^{\beta}  \pi_k^\gamma 
- \delta_k^i \delta_\alpha^{\gamma} \pi_j^\beta\Bigr).
 \end{array} 
\renewcommand{\arraystretch}{1}
\end{equation}

We will eliminate the fiber forms 
$\pi_{\alpha k}^{\beta \gamma}, \pi_{jk}^{i\gamma}$, and $\pi_k^\gamma$ from equations (3.34). 
To this end, we construct the following three objects:
\begin{equation}\label{eq:3.35}
\renewcommand{\arraystretch}{1.5}
 \begin{array}{ll}
x^{i\beta\gamma}_{\alpha jk} 
=& \!\!\!\!  - \displaystyle \frac{2}{q^2 - 1} 
\delta_{[j}^{i}\Bigl(pu_{|\alpha| k]}^{[\beta\gamma]}+
 u_{|\alpha|k]}^{[\gamma\beta]}\Bigr),  \\
y^{i\beta\gamma}_{\alpha jk} 
=& \!\!\!\! - \displaystyle \frac{2}{p^2 - 1} 
 \delta_{\alpha}^{[\beta}
\Bigl(qu_{[jk]}^{|i|\gamma]} + 
 u_{[kj]}^{|i|\gamma]}\Bigr), \\ 
z^{i\beta\gamma}_{\alpha jk} 
=& \!\!\!\! \displaystyle \frac{2}{(p^{2}-1)(q^{2}-1)}
 \Bigl[(pq-1)\Bigl(\delta_{[j}^{i}\delta_{|\alpha|}^{[\beta} 
 u_{|l|k]}^{|l|\gamma]}
 + \delta_{[k}^{i}\delta_{|\alpha|}^{[\beta} 
 u_{j]l}^{|l|\gamma]}\Bigr)  \\
& \!\!\!\! + (p-q)\Bigl(\delta_{[j}^{i}\delta_{|\alpha|}^{[\beta} 
 u_{k]l}^{|l|\gamma]} 
+ \delta_{[k}^{i} \delta_{|\alpha|}^{[\beta} 
u_{|l|j]}^{|l|\gamma]}\Bigr) \Bigr],
 \end{array} 
\renewcommand{\arraystretch}{1}
\end{equation}
where the quantities $u^{\beta\gamma}_{\alpha k}$ and 
$u^{i\gamma}_{jk}$ are defined by formulas (3.21). 
A straightforward calculation with help of equations (3.22) 
and (3.23) gives the following differential  equations 
for the objects $x^{i\beta\gamma}_{\alpha jk},\; y^{i\beta\gamma}_{\alpha jk}$, and $z^{i\beta\gamma}_{\alpha jk}$:
\begin{equation}\label{eq:3.36}
\renewcommand{\arraystretch}{1.5}
 \begin{array}{ll}
\nabla_\delta x^{i\beta\gamma}_{\alpha jk} 
+ x^{i\beta\gamma}_{\alpha jk} \pi 
=&  \!\!\!\! \displaystyle \frac{1}{2} \Bigl[
\delta_j^i (\pi_{\alpha k}^{\beta \gamma}  
- \delta_\alpha^\beta  \pi_k^\gamma) 
- \delta_k^i (\pi_{\alpha j}^{\gamma \beta} 
- \delta_\alpha^{\gamma} \pi_j^\beta)\Bigr] \\
&\!\!\!\! + \displaystyle \frac{1}{2 (q^2 - 1)} \Bigl[
q (\delta_k^i \delta_\alpha^\beta \pi_{jl}^{l\gamma}
- \delta_j^i \delta_\alpha^{\gamma}   \pi_{kl}^{l\beta})
 + \delta_k^i \delta_\alpha^{\gamma} \pi_{jl}^{l\beta} 
-  \delta_j^i \delta_\alpha^\beta  \pi_{kl}^{l\gamma}\Bigr], \\ 
\nabla_\delta y^{i\beta\gamma}_{\alpha jk} 
+ y^{i\beta\gamma}_{\alpha jk} \pi 
= & \!\!\!\! \displaystyle \frac{1}{2} \Bigl[
\delta_\alpha^\beta (\pi_{j k}^{i \gamma}  
- \delta_j^i  \pi_k^\gamma) 
- \delta_\alpha^\gamma (\pi_{kj}^{i\beta} 
- \delta_k^i \pi_j^\beta \Bigr] \\
&\!\!\!\! + \displaystyle \frac{1}{2 (p^2 - 1)} \Bigl[
p (\delta_j^i \delta_\alpha^\gamma \pi_{\delta k}^{\beta \delta} 
- \delta_k^i \delta_\alpha^\beta  \pi_{\delta j}^{\gamma \delta})
 + \delta_k^i \delta_\alpha^\gamma \pi_{\delta j}^{\beta \delta} 
-  \delta_j^i \delta_\alpha^\beta 
\pi_{\delta k}^{\gamma \delta}\Bigr], \\ 
\nabla_\delta z^{i\beta\gamma}_{\alpha jk} 
+ z^{i\beta\gamma}_{\alpha jk} \pi 
=& \!\!\!\! - \displaystyle \frac{1}{2(p^{2}-1)} 
\Bigl[p (\delta_j^i \delta_\alpha^\gamma 
\pi_{\delta k}^{\beta \delta} 
- \delta_k^i \delta_\alpha^\beta  \pi_{\delta j}^{\gamma \delta})
 + \delta_k^i \delta_\alpha^\gamma \pi_{\delta j}^{\beta \delta} 
-  \delta_j^i \delta_\alpha^\beta 
\pi_{\delta k}^{\gamma \delta}\Bigr] \\
&\!\!\!\! - \displaystyle \frac{1}{2(q^{2}-1)} 
\Bigl[q (\delta_k^i \delta_\alpha^\beta \pi_{jl}^{l\gamma}
- \delta_j^i \delta_\alpha^{\gamma}   \pi_{kl}^{l\beta})
 + \delta_k^i \delta_\alpha^{\gamma} \pi_{jl}^{l\beta} 
- \delta_j^i \delta_\alpha^\beta  
\pi_{kl}^{l\gamma}\Bigr]\\
&\!\!\!\!  - \displaystyle \frac{1}{2} 
\Bigl(\delta_k^i \delta_\alpha^{\gamma} \pi_j^\beta  
-  \delta_j^i \delta_\alpha^\beta \pi_k^\gamma\Bigr). 
 \end{array} 
\renewcommand{\arraystretch}{1}
\end{equation}

If we set 
\begin{equation}\label{eq:3.37}
\widetilde{a}^{i\beta\gamma}_{\alpha jk} 
= u^{i\beta\gamma}_{\alpha jk} 
+  x^{i\beta\gamma}_{\alpha jk}
+  y^{i\beta\gamma}_{\alpha jk}
+ z^{i\beta\gamma}_{\alpha jk},
\end{equation}
then by (3.34) and (3.36), it is easy to check that 
$$
\nabla_\delta \widetilde{a}^{i\beta\gamma}_{\alpha jk} 
+ \widetilde{a}^{i\beta\gamma}_{\alpha jk} \pi = 0.
$$
This means that the quantities 
$\widetilde{a}^{i\beta\gamma}_{\alpha jk}$ form a relative tensor 
of weight $-1$. 
Using (3.12), it is easy to verify that 
the tensor $\widetilde{a}_{\alpha jk}^{i\beta\gamma}$ satisfies 
the conditions similar to conditions (3.30) and (3.31). 

We will prove now the following proposition: 

\begin{proposition} The relative tensor 
$\widetilde{a} = \{\widetilde{a}_{\alpha jk}^{i\beta\gamma}\}$ 
defined by formulas $(3.37)$ and $(3.35)$ and satisfying 
the conditions similar to conditions $(3.30)$ and $(3.31)$ 
coincides with the tensor $a = \{a_{\alpha jk}^{i\beta\gamma}\}$: 
$\widetilde{a} = a$. 
\end{proposition}

{\sf Proof.} Let us assume that
 from the beginning we have the tensor  \linebreak 
 $a = \{a_{\alpha jk}^{i\beta\gamma}\}$ in place of the 
object 
$u_{\alpha jk}^{i\beta\gamma}$ and that tensor satisfies 
 conditions (3.30) and (3.31); in other words, we have 
equations 
(3.30)--(3.33) and (3.29). Then, applying  
(3.30), we find from  (3.35) that 
$x_{\alpha jk}^{i\beta\gamma} 
= y_{\alpha jk}^{i\beta\gamma} 
= z_{\alpha jk}^{i\beta\gamma} = 0$, and consequently 
from formula (3.37) we find that 
$\widetilde{a}_{\alpha jk}^{i\beta\gamma} 
= a_{\alpha jk}^{i\beta\gamma}$. \rule{3mm}{3mm}

By Proposition 3.3, 
if we substitute for $x^{i\beta\gamma}_{\alpha jk}, 
y^{i\beta\gamma}_{\alpha jk}$ and $z^{i\beta\gamma}_{\alpha jk}$ 
in equations (3.37) their values  (3.35), we find 
that 
\begin{equation}\label{eq:3.38}
\renewcommand{\arraystretch}{1.5}
\begin{array}{lll}
a_{\alpha jk}^{i\beta\gamma} & \!\!\!\!  = \!\!\!\!  & 
u_{\alpha jk}^{i\beta\gamma}-
\displaystyle \frac{2}{q^{2}-1} 
\delta_{[j}^{i}\Bigl(qu_{|\alpha|k]}^{[\beta\gamma]}+
 u_{|\alpha|k]}^{[\gamma\beta]}\Bigr)  
  -\displaystyle \frac{2}{p^{2}-1}\delta_{\alpha}^{[\beta}
\Bigl(pu_{[jk]}^{|i|\gamma]}+
 u_{[kj]}^{|i|\gamma]}\Bigr)  \\
      &\!\!\!\! \!\!\!\!     & 
+ \displaystyle \frac{2}{(p^{2}-1)(q^{2}-1)}
 \Bigl[(pq-1)\Bigl(\delta_{[j}^{i}\delta_{|\alpha|}^{[\beta}
 u_{k]}^{\gamma]} 
+ \delta_{[k}^i\delta_{|\alpha|}^{[\beta}
 \widetilde{u}_{j]}^{\gamma]}\Bigr)  \\
      &\!\!\!\! \!\!\!\!     & 
+ (q-p)\Bigl(\delta_{[j}^i\delta_{|\alpha|}^{[\beta} 
 \widetilde{u}_{k]}^{\gamma]}
+ \delta_{[k}^i\delta_{|\alpha|}^{[\beta} 
u_{j]}^{\gamma]}\Bigr) \Bigr],
\end{array}
\renewcommand{\arraystretch}{1}
\end{equation}
where the alternation is carried out with respect 
to the pairs of indices ${\beta \choose j}$, ${\gamma \choose k}$ 
or ${\beta \choose k}$, ${\gamma \choose j}$, and 
$
u_k^\gamma = u_{lk}^{l\gamma} = u_{\sigma lk}^{l\sigma\gamma} 
= u_{\sigma k}^{\sigma\gamma}, \;\;
\widetilde{u}_k^\gamma = u_{kl}^{l\gamma} 
= u_{\sigma kl}^{l\sigma\gamma} = - u_{\sigma lk}^{l\gamma\sigma} 
\linebreak = - u_{\sigma k}^{\gamma\sigma}.
$

The expression (3.38) of the components of the tensor 
$a = \{a_{\alpha jk}^{i\beta\gamma}\}$ in the general (not 
reduced) third-order frame was found  by Goldberg in 
~\cite{G75} (see also \cite{G88}, \S 2.2). Note that using another method, 
Hangan \cite{H80} deduced this expression again. 

\section{The Complete Structure Object 
of an Almost \newline Grassmann Manifold}

\setcounter{equation}{0}

{\bf 1.} In this section we will construct  the second 
and the third structural objects of the almost Grassmann 
structure.

 From equations (3.29) it follows that, after 
the reduction of the second-order frame bundle made 
in Section {\bf 3},   we have the following equations on the 
manifold $AG (p - 1, p + q - 1)$:
\begin{equation}\label{eq:4.1}
\renewcommand{\arraystretch}{1.5}
 \begin{array}{ll}
\omega_{\alpha k}^{\beta \gamma} 
= \displaystyle \frac{q}{p+q} (\delta_\alpha^\beta 
 \omega_{k}^{\gamma} - p \delta_\alpha^\gamma \omega_{k}^{\beta}) + \widehat{w}^{\beta\gamma\delta}_{\alpha kl} \omega_\delta^l, \\ 
\omega_{jk}^{i \gamma} 
= \displaystyle \frac{p}{p+q} (\delta_j^i 
 \omega_k^\gamma - q \delta_k^i \pi_j^\gamma) 
+ \widehat{w}^{i\gamma\delta}_{jkl} \omega_\delta^l. 
\end{array} 
\renewcommand{\arraystretch}{1}
\end{equation}
Substituting these forms into  system (3.15), we obtain
\begin{equation}\label{eq:4.2}
\renewcommand{\arraystretch}{1.5}
 \begin{array}{ll}
d \omega_\alpha^\beta - \omega_\alpha^\gamma \wedge \omega_\gamma^\beta 
= \displaystyle \frac{q}{p+q}(\delta_\alpha^\beta \omega_\gamma^k 
\wedge \omega^\gamma_k - p \omega_\alpha^k \wedge \omega_k^\beta) 
+  w_{\alpha kl}^{\beta\gamma\delta}  \omega_\gamma^k \wedge 
\omega_\delta^l, \\
 d \omega_j^i - \omega_j^k  \wedge \omega_k^i
= \displaystyle \frac{p}{p+q}(\delta_j^i \omega^\gamma_k 
\wedge \omega_\gamma^k 
- q \omega_j^\gamma  \wedge \omega_\gamma^i) 
+  w_{jkl}^{i\gamma\delta}  \omega_\gamma^k \wedge 
\omega_\delta^l, \\
d \omega = \omega_i^\alpha \wedge \omega_\alpha^i,
\end{array} 
\renewcommand{\arraystretch}{1}
\end{equation}
where $w_{\alpha kl}^{\beta\gamma\delta} 
\! =\!  \widehat{w}_{\alpha [kl]}^{\beta[\gamma\delta]}$ 
and $w_{jkl}^{i\gamma\delta} \! = \! \widehat{w}_{j[kl]}^{i[\gamma\delta]}$. 
The latter quantities satisfy the relations
\begin{equation}\label{eq:4.3}
  w_{\alpha kl}^{\beta\gamma\delta} 
= -   w_{\alpha lk}^{\beta\delta\gamma}, \;\; 
  w_{jkl}^{i\gamma\delta} = - w_{jlk}^{i\delta\gamma}
\end{equation}
and also the relations
\begin{equation}\label{eq:4.4}
  w_{\alpha kl}^{\alpha\gamma\delta} = 0,  \;\; 
  w_{ikl}^{i\gamma\delta} = 0,
\end{equation}
which follow from equations (4.1) and (3.17).

 Exterior differentiation of equations (4.2) 
gives the following exterior cubic equations:
\begin{equation}\label{eq:4.5}
\renewcommand{\arraystretch}{1.5}
 \begin{array}{ll}
(\nabla w^{\beta\gamma\delta}_{\alpha kl} 
+ 2 w^{\beta\gamma\delta}_{\alpha kl} \omega) 
\wedge \omega_\gamma^k  \wedge \omega^l_\delta 
- 2 w^{\beta\gamma\delta}_{\alpha kl} 
a^{l \varepsilon\zeta}_{\delta mn}   \omega_\gamma^k
 \wedge \omega_\varepsilon^m  \wedge \omega_\zeta^n \\
+ \displaystyle \frac{pq}{p+q}  
\biggl[(\Omega^\beta_k + \omega \wedge \omega_k^\beta) 
\wedge \omega_\alpha^k - 
a^{m\gamma\delta}_{\alpha kl} \omega_m^\beta  
\wedge \omega_\gamma^k  \wedge \omega^l_\delta\biggr]
 = 0,
\end{array} 
\renewcommand{\arraystretch}{1}
\end{equation}
\begin{equation}\label{eq:4.6}
\renewcommand{\arraystretch}{1.5}
 \begin{array}{ll}
(\nabla w^{i\gamma\delta}_{jkl} 
+ 2 w^{i\gamma\delta}_{jkl} \omega)  
\wedge \omega_\gamma^k  \wedge \omega^l_\delta 
- 2 w^{i\gamma\delta}_{jkl} 
a^{l \varepsilon\zeta}_{\delta mn}   \omega_\gamma^k
 \wedge \omega_\varepsilon^m  \wedge \omega_\zeta^n \\
- \displaystyle \frac{pq}{p+q} 
\biggl[(\Omega^\gamma_j + \omega \wedge \omega_j^\gamma) \wedge 
 \omega_\gamma^i  
- a^{i\delta\varepsilon}_{\gamma lm} \omega_j^\gamma 
\wedge \omega_\delta^l  \wedge \omega^m_\varepsilon\biggr]  = 0,
\end{array} 
\renewcommand{\arraystretch}{1}
\end{equation}
and 
\begin{equation}\label{eq:4.7}
(\Omega_i^\alpha + \omega \wedge \omega_i^\alpha) \wedge 
 \omega_\alpha^i  - a^{i\beta\gamma}_{\alpha jk} 
 \omega^\alpha_i  \wedge \omega_\beta^j  \wedge \omega_\gamma^k 
= 0,
\end{equation}
where 
$$
\renewcommand{\arraystretch}{1.5}
 \begin{array}{ll}
\Omega_i^\alpha = d \omega_i^\alpha - \omega_i^\gamma 
\wedge \omega_\gamma^\alpha - \omega_i^l \wedge 
\omega_l^\alpha,\\
\nabla w^{\beta\gamma\delta}_{\alpha kl} 
= d w^{\beta\gamma\delta}_{\alpha kl} 
- w^{\beta\gamma\delta}_{\varepsilon kl} \omega^\varepsilon_\alpha 
- w^{\beta\gamma\delta}_{\alpha jl} \omega^j_k 
- w^{\beta\gamma\delta}_{\alpha kj} \omega^j_l 
+ w^{\varepsilon\gamma\delta}_{\alpha kj} \omega^\beta_\varepsilon 
+ w^{\beta\varepsilon\delta}_{\alpha kj} \omega^\gamma_\varepsilon  
+ w^{\beta\gamma\varepsilon}_{\alpha kj} \omega^\delta_\varepsilon, \\  
\nabla w^{i\gamma\delta}_{jkl} 
= d w^{i\gamma\delta}_{jkl} 
- w^{i\gamma\delta}_{mkl}  \omega^m_j 
- w^{i\gamma\delta}_{jml}  \omega^m_k  
- w^{i\gamma\delta}_{jkm} \omega^m_l 
+ w^{m\gamma\delta}_{jkl}  \omega_m^i 
+ w^{i\varepsilon\delta}_{jkl}  \omega^\varepsilon_\gamma 
+ w^{i\gamma\varepsilon}_{jkl}  \omega^\varepsilon_\delta. 
\end{array} 
\renewcommand{\arraystretch}{1}
$$

 From equations (4.5)--(4.7) it follows that the 2-form 
$\Omega_i^\alpha + \omega \wedge \omega_i^\alpha$ can be 
expressed as follows:
\begin{equation}\label{eq:4.8}
\Omega_i^\alpha + \omega \wedge \omega_i^\alpha
=  \omega^{\alpha\beta}_{ij}  \wedge \omega_\beta^j.
\end{equation}
The forms $\omega^{\alpha\beta}_{ij}$ are defined on the 
fourth-order frame bundle associated with the almost Grassmann 
structure $AG (p - 1, p + q - 1)$.

Substituting expressions (4.8) into (4.5)--(4.7), 
we obtain 
\begin{equation}\label{eq:4.9}
\renewcommand{\arraystretch}{1.5}
 \begin{array}{ll}
\Bigl[\nabla w^{\beta\gamma\delta}_{\alpha kl} 
&+ 2 w^{\beta\gamma\delta}_{\alpha kl} \omega 
- \displaystyle \frac{pq}{p+q} \bigl(\delta_\alpha^\gamma 
\omega^{\beta\delta}_{kl} + 
a^{m\gamma\delta}_{\alpha kl} \omega_m^\beta\bigr)\Bigr] 
\wedge \omega_\gamma^k  \wedge \omega^l_\delta \\
&- 2 w^{\beta\varepsilon\sigma}_{\alpha ms} 
a^{s \gamma\delta}_{\sigma kl}  \wedge \omega_\varepsilon^m
 \wedge \omega_\gamma^k  \wedge \omega_\delta^l = 0,
\end{array} 
\renewcommand{\arraystretch}{1}
\end{equation}
\begin{equation}\label{eq:4.10}
\renewcommand{\arraystretch}{1.5}
 \begin{array}{ll}
\Bigl[\nabla w^{i\gamma\delta}_{jkl} 
&+ 2 w^{i\gamma\delta}_{jkl} \omega 
+ \displaystyle \frac{pq}{p+q} 
\bigl(\delta_k^i \omega_{jl}^{\gamma\delta} + 
a^{i\gamma\delta}_{\varepsilon kl} 
\omega_j^\varepsilon\bigr)\Bigr] 
\wedge \omega_\gamma^k  \wedge \omega^l_\delta \\
&- 2 w^{i\varepsilon\sigma}_{jms} 
a^{s \gamma\delta}_{\sigma kl}  \wedge \omega_\varepsilon^m
 \wedge \omega_\gamma^k  \wedge \omega_\delta^l = 0,
\end{array} 
\renewcommand{\arraystretch}{1}
\end{equation}
and 
\begin{equation}\label{eq:4.11}
(\omega_{ij}^{\alpha\beta}  
+ a^{k\alpha\beta}_{\gamma ij} \omega_k^\gamma) 
\wedge \omega_\alpha^i  \wedge \omega_\beta^j = 0.
\end{equation}
It follows from (4.9)--(4.11) that for 
$\omega_\alpha^i = 0$, we have 
\begin{equation}\label{eq:4.12}
\nabla_\delta w^{\beta\gamma\delta}_{\alpha kl} 
+ 2 w^{\beta\gamma\delta}_{\alpha kl} \pi 
- \displaystyle \frac{pq}{p+q} \bigl(\delta_\alpha^{[\gamma} 
\pi^{|\beta|\delta]}_{[kl]} +  
a^{m\gamma\delta}_{\alpha kl} \pi_m^\beta\bigr) = 0,
\end{equation}
\begin{equation}\label{eq:4.13}
\nabla_\delta w^{i\gamma\delta}_{jkl} 
+ 2 w^{i\gamma\delta}_{jkl} \pi 
+ \displaystyle \frac{pq}{p+q} 
\bigl(\delta_{[k}^i \pi_{|j|l]}^{[\gamma\delta]} + 
a^{i\gamma\delta}_{\varepsilon kl} \pi_j^\varepsilon\bigr) = 0,
\end{equation}
and 
\begin{equation}\label{eq:4.14}
\pi_{[kl]}^{[\gamma\delta]}  
+ a^{m\gamma\delta}_{\varepsilon kl} \pi_m^\varepsilon = 0,
\end{equation}
where the alternation is carried over the pairs of indices 
$
{\gamma \choose k}$ and 
$
{\delta \choose l}$ and 
as usual $\pi_k^\gamma = \omega_k^\gamma (\delta), \;
\pi_{jk}^{\beta\gamma} = \omega_{jk}^{\beta\gamma} (\delta)$, 
and $\pi = \omega (\delta)$.

 From (4.14) it follows that the form 
$\pi_{kl}^{\gamma\delta}$ can be written as 
\begin{equation}\label{eq:4.15}
\pi_{kl}^{\gamma\delta}  = \widetilde{\pi}_{kl}^{\gamma\delta} 
-  a^{m\gamma\delta}_{\varepsilon kl} \pi_m^\varepsilon,
\end{equation}
where 
\begin{equation}\label{eq:4.16}
\widetilde{\pi}_{kl}^{\gamma\delta} 
= \widetilde{\pi}_{lk}^{\delta\gamma}. 
\end{equation}
The 1-forms $\widetilde{\pi}_{kl}^{\gamma\delta}$ determine 
admissible transformations of third order frames associated with 
the almost Grassmann structure $AG (p - 1, p + q - 1)$.

If we apply equation (4.15), we can write equations 
(4.12) and (4.13) in the form:
\begin{equation}\label{eq:4.17}
\renewcommand{\arraystretch}{1.5}
 \begin{array}{ll}
\nabla_\delta w^{\beta\gamma\delta}_{\alpha kl} 
+ 2 w^{\beta\gamma\delta}_{\alpha kl} \pi 
&- \displaystyle \frac{pq}{2(p+q)} 
\Bigl[\delta_\alpha^{\gamma} 
\widetilde{\pi}^{\beta\delta}_{kl} 
- \delta_\alpha^{\delta} 
\widetilde{\pi}^{\beta\gamma}_{lk} \\
&+ (2\delta_\varepsilon^\beta a^{m\gamma\delta}_{\alpha kl} 
- \delta_\alpha^\gamma a^{m\beta\delta}_{\varepsilon kl} 
+ \delta_\alpha^\delta a^{m\beta\gamma}_{\varepsilon lk})  
\pi_m^\varepsilon\Bigr] = 0,
\end{array} 
\renewcommand{\arraystretch}{1}
\end{equation}
\begin{equation}\label{eq:4.18}
\renewcommand{\arraystretch}{1.5}
 \begin{array}{ll}
\nabla_\delta w^{i\gamma\delta}_{jkl} 
+ 2 w^{i\gamma\delta}_{jkl} \pi 
&+ \displaystyle \frac{pq}{2(p+q)} 
\Bigl[\delta_k^i \widetilde{\pi}_{jl}^{\gamma\delta} 
- \delta_l^i \widetilde{\pi}_{jk}^{\delta\gamma} \\
& + (2\delta_j^m a^{i\gamma\delta}_{\varepsilon kl}
- \delta_k^i a^{m\gamma\delta}_{\varepsilon jl} 
+ \delta_l^i a^{m\delta\gamma}_{\varepsilon jk})
 \pi_m^\varepsilon\Bigr] = 0.
\end{array} 
\renewcommand{\arraystretch}{1}
\end{equation}

{\bf 2.} 
If we contract equation (4.17) with respect to the indices 
$\alpha$ and $\gamma$,  equation (4.18) with respect to the 
indices $i$ and $k$ and change the notation of some indices, we 
obtain
\begin{equation}\label{eq:4.19}
\nabla_\delta w^{\gamma\alpha\delta}_{\alpha kl} 
+ 2 w^{\gamma\alpha\delta}_{\alpha kl} \pi 
- \displaystyle \frac{pq}{2(p+q)} 
\Bigl[p \widetilde{\pi}^{\gamma\delta}_{kl} 
- \widetilde{\pi}^{\delta\gamma}_{kl} 
+ (a^{m\gamma\delta}_{\varepsilon lk}
- p a^{m\gamma\delta}_{\varepsilon kl}) 
\pi_m^\varepsilon\Bigr] = 0,
\end{equation}
\begin{equation}\label{eq:4.20}
\nabla_\delta w^{i\gamma\delta}_{kil} 
+ 2 w^{i\gamma\delta}_{kil} \pi 
+ \displaystyle \frac{pq}{2(p+q)} 
\Bigl[q \widetilde{\pi}_{kl}^{\gamma\delta} 
-  \widetilde{\pi}_{kl}^{\delta\gamma} 
 + (a^{m\delta\gamma}_{\varepsilon kl} 
- q a^{m\gamma\delta}_{\varepsilon kl}) 
  \pi_m^\varepsilon\Bigr] = 0.
\end{equation}

Using the quantities $w^{\gamma\alpha\delta}_{\alpha kl}$ 
and $w^{i\gamma\delta}_{kil}$, we now construct the following 
new object:
\begin{equation}\label{eq:4.21}
 w^{\gamma\delta}_{kl} 
= w^{\gamma\alpha\delta}_{\alpha kl} 
-  w^{i\gamma\delta}_{kil}  
+ w^{\delta\alpha\gamma}_{\alpha lk} 
- w^{i\delta\gamma}_{lik}.
\end{equation}
By means of (4.19), (4.20), and (4.16), 
it is easy to prove that the quantities 
$w^{\gamma\delta}_{kl}$ defined by (4.21) satisfy the following 
differential equations:
\begin{equation}\label{eq:4.22}
\nabla_\delta w^{\gamma\delta}_{kl} + 2 w^{\gamma\delta}_{kl} \pi 
= \displaystyle \frac{pq}{p+q} 
\bigl[(p+q) \widetilde{\pi}_{kl}^{\gamma\delta} 
- 2 \widetilde{\pi}_{kl}^{\delta\gamma}\bigr]. 
\end{equation}
Formulas (4.22) show that  the quantities 
$w^{\gamma\delta}_{kl}$ form a geometric object which is defined 
in a third-order differential neighborhood of 
the almost Grassmann structure $AG (p - 1, p + q - 1)$. 

Let us prove that if we make 
a specialization of fourth order frames, then we can reduce 
the geometric object $w^{\gamma\delta}_{kl}$ defined by 
(4.21) to 0. 
In our proof we will again apply the  method  used in Section 
 {\bf 3}.

 We must show that the 
1-forms on the right-hand sides of equations (4.22) are 
linearly independent. First, we note that by (4.16) and 
(4.21), both the components 
$w^{\gamma\delta}_{kl}$ and the forms $\widetilde{\pi}_{kl}^{\gamma\delta}$ are symmetric with respect 
to the pairs of indices $
{\gamma \choose k}$ and 
$
{\delta \choose l}$, and there are 
$
\frac{1}{2} pq(pq+1)$ 
linearly independent forms among the 1-forms 
$\widetilde{\pi}_{kl}^{\gamma\delta}$. Let us equate 
to 0 the right-hand sides of equations (4.22):
\begin{equation}\label{eq:4.23}
(p+q) \widetilde{\pi}_{kl}^{\gamma\delta} 
- 2 \widetilde{\pi}_{kl}^{\delta\gamma} = 0.
\end{equation}
Interchanging  the indices  $\gamma$ and 
$\delta$ in (4.23), we obtain
\begin{equation}\label{eq:4.24}
- 2 \widetilde{\pi}_{kl}^{\gamma\delta}
+ (p+q) \widetilde{\pi}_{kl}^{\delta\gamma} = 0.
\end{equation}
Since the determinant of the system of equations 
(4.23)--(4.24) is equal to $(p+q)^2 - 4 \neq 0$, the 
system has only the trivial solution.

There are $
\frac{1}{2} pq(pq+1)$ 
linearly independent forms among the 1-forms  
 $(p+q) \widetilde{\pi}_{kl}^{\gamma\delta} 
- 2 \widetilde{\pi}_{kl}^{\delta\gamma}$. But, up to the factor 
$
\frac{pq}{p+q}$, the latter forms 
 coincide with the right-hand sides of equations (4.22). 

Hence the geometric object $w^{\gamma\delta}_{kl}$ 
can be reduced to 0. This operation leads to a reduction 
of the set of third-order frames of 
the almost Grassmann structure $AG (p - 1, p + q - 1)$. 
Before this reduction, the set of third-order frames 
depended on $\frac{1}{2} pq(pq+1)$ 
parameters equal 
to the number of linearly independent forms among the forms 
$\widetilde{\pi}_{kl}^{\gamma\delta}$. 
After the reduction, the forms 
$\widetilde{\pi}_{kl}^{\gamma\delta}$  vanish, 
\begin{equation}\label{eq:4.25}
 \widetilde{\pi}_{kl}^{\gamma\delta} = 0,
\end{equation} 
and the forms $\pi_{kl}^{\gamma\delta}$  are expressed 
in terms of the 1-forms $\pi_k^\gamma$, 
\begin{equation}\label{eq:4.26}
\pi_{kl}^{\gamma\delta}  
= -  a^{m\gamma\delta}_{\varepsilon kl} \pi_m^\varepsilon.
\end{equation}
But as we saw earlier, the forms $\pi_k^\gamma$ determine 
admissible transformations of second-order frames. 
This means that after the above reduction
the group of admissible transformations of 
third-order frames coincides with the group of admissible 
transformations of second-order frames. This implies the 
following result:

\begin{theorem}
The almost  Grassmann manifold 
$AG (p-1, p+q-1)$ is a $G$-structure of finite type two. 
\end{theorem}
This result is analogous to the result for the conformal 
$CO (p, q)$-structures (see \cite{AG96},  \S 4.1).  The theorem 
similar to Theorem 4.1 was proved in \cite{H80} in terms of Lie 
algebras (see \cite{St64}, Ch. 7, \S 3, for the definition of 
a $G$-structure of finite type).

Denote  the values of the 
 quantities $w^{\beta\gamma\delta}_{\alpha kl}$ 
and $w^{i\gamma\delta}_{jkl}$ in the reduced fourth-order frames 
by $b^{\beta\gamma\delta}_{\alpha kl}$ 
and $b^{i\gamma\delta}_{jkl}$, respectively. 
Then the quantities $b^{\beta\gamma\delta}_{\alpha kl}$ 
and $b^{i\gamma\delta}_{jkl}$ satisfy 
the differential equations  obtained from 
equations (4.17) and (4.18) by means of (4.25):

\begin{equation}\label{eq:4.27}
\nabla_\delta b^{\beta\gamma\delta}_{\alpha kl} 
+ 2 b^{\beta\gamma\delta}_{\alpha kl} \pi 
- \displaystyle \frac{pq}{2(p+q)} 
\Bigl(2\delta_\varepsilon^\beta a^{m\gamma\delta}_{\alpha kl} 
- \delta_\alpha^\gamma a^{m\beta\delta}_{\varepsilon kl} 
+ \delta_\alpha^\delta a^{m\beta\gamma}_{\varepsilon lk}\Bigr)  
\pi_m^\varepsilon = 0 
\end{equation}
and 
\begin{equation}\label{eq:4.28}
\nabla_\delta b^{i\gamma\delta}_{jkl} 
+ 2 b^{i\gamma\delta}_{jkl} \pi 
+ \displaystyle \frac{pq}{2(p+q)} 
\Bigl(2\delta_j^m a^{i\gamma\delta}_{\varepsilon kl}
- \delta_k^i a^{m\gamma\delta}_{\varepsilon jl} 
+ \delta_l^i a^{i\delta\gamma}_{\varepsilon jk}\Bigl)
 \pi_m^\varepsilon = 0.
\end{equation}
They also satisfy  the conditions (cf. (4.21)) 
\begin{equation}\label{eq:4.29}
 b^{\gamma\alpha\delta}_{\alpha kl} 
-  b^{i\gamma\delta}_{kil}  
+ b^{\delta\alpha\gamma}_{\alpha lk} 
- b^{i\delta\gamma}_{lik} = 0
\end{equation}
and
\begin{equation}\label{eq:4.30}
  b_{\alpha kl}^{\beta\gamma\delta} 
= - b_{\alpha lk}^{\beta\delta\gamma}, \;\; 
  b_{jkl}^{i\gamma\delta} = - b_{jlk}^{i\delta\gamma},
\end{equation}
and  the relations
\begin{equation}\label{eq:4.31}
  b_{\alpha kl}^{\alpha\gamma\delta} = 0,  \;\; 
  b_{ikl}^{i\gamma\delta} = 0.
\end{equation}
The relations (4.30) and (4.31) follow from condition 
(4.3) and (4.4). 

Equations (4.26) show that the 1-forms 
$\omega_{kl}^{\gamma\delta}$ occurring in equations 
(4.8) are expressed in terms of the forms $\omega_k^\gamma$ 
and the basis forms $\omega_\gamma^k$: 
\begin{equation}\label{eq:4.32}
\omega_{ij}^{\alpha\beta}  
= -  a^{k\alpha\beta}_{\gamma ij} \omega_k^\gamma 
+ \widehat{c}^{\alpha\beta\gamma}_{ijk} \omega_\gamma^k.
\end{equation}
This means that the fiber forms $\omega_{ij}^{\alpha\beta}$
 associated with the fourth-order 
frame bundle are expressed in terms of the fiber forms 
$\omega_k^\gamma$, defined on the third-order frame bundle, and 
the basis forms $\omega_\gamma^k$. 

{\bf 3.} Substituting for  the forms 
$\omega_{ij}^{\alpha\beta}$ in equations (4.8) 
their values  (4.32), we find that
\begin{equation}\label{eq:4.33}
d \omega_i^\alpha - \omega_i^\beta \wedge \omega_\beta^\alpha 
- \omega_i^j \wedge \omega_j^\alpha  
+ \omega \wedge \omega_i^\alpha
=    c_{ijk}^{\alpha\beta\gamma} 
 \omega^k_\gamma \wedge \omega_\beta^j 
- a_{\gamma ij}^{k\alpha\beta} \omega_k^\gamma \wedge 
\omega_\beta^j,
\end{equation}
where $c_{ijk}^{\alpha\beta\gamma} 
= \widehat{c}_{i[jk]}^{\alpha[\beta\gamma]}$, and 
 the alternation is carried over  the vertical pairs of indices.

If we substitute for the forms $\omega_{ij}^{\alpha\beta}$ 
in equations (4.11) their values  (4.32), we find 
that  the quantities $c_{ijk}^{\alpha\beta\gamma}$ 
must satisfy the following condition:
\begin{equation}\label{eq:4.34}
  c_{[ijk]}^{[\alpha\beta\gamma]} = 0.
\end{equation}
Equations (4.33) together with equations (3.33) and 
(4.2) make up the complete system of the 
structure equations of the almost Grassmann structure 
$AG (p - 1, p + q - 1)$:
\begin{equation}\label{eq:4.35}
\renewcommand{\arraystretch}{1.5}
\left\{
 \begin{array}{ll}
 d\omega_\alpha^i - \omega_\alpha^j \wedge 
\omega_j^i - \omega_\alpha^\beta \wedge \omega_\beta^i 
-  \omega \wedge \omega_\alpha^i = a_{\alpha jk}^{i\beta\gamma} 
\omega_{\beta}^{j} \wedge \omega_{\gamma}^{k}, \\
d \omega_\alpha^\beta - \omega_\alpha^\gamma \wedge \omega_\gamma^\beta 
= \displaystyle \frac{q}{p+q}\Bigl(\delta_\alpha^\beta 
\omega_\gamma^k \wedge \omega^\gamma_k - p \omega_\alpha^k 
\wedge \omega_k^\beta\Bigr) 
+  b_{\alpha kl}^{\beta\gamma\delta}  \omega_\gamma^k \wedge 
\omega_\delta^l, \\
 d \omega_j^i - \omega_j^k  \wedge \omega_k^i
= \displaystyle \frac{p}{p+q} \Bigl(\delta_j^i \omega^\gamma_k 
\wedge \omega_\gamma^k 
- q \omega_j^\gamma  \wedge \omega_\gamma^i\Bigr) 
+  b_{jkl}^{i\gamma\delta}  \omega_\gamma^k \wedge 
\omega_\delta^l, \\
d \omega = \omega_i^\alpha \wedge \omega_\alpha^i,\\
d \omega_i^\alpha - \omega_i^\beta \wedge \omega_\beta^\alpha 
- \omega_i^j \wedge \omega_j^\alpha  
+ \omega \wedge \omega_i^\alpha
=    c_{ijk}^{\alpha\beta\gamma} 
 \omega^k_\gamma \wedge \omega_\beta^j 
- a_{\gamma ij}^{k\alpha\beta} \omega_k^\gamma \wedge 
\omega_\beta^j.
 \end{array}
\right.
\renewcommand{\arraystretch}{1}
\end{equation}

We have proved above that a structure 
$AG (p-1, p+q-1)$ is a $G$-structure of finite type two. The 
invariant forms of this structure are divided into three groups: 
$\{\omega_\alpha^i\}, \{\omega^i_j, 
\omega_\alpha^\beta, \omega\}$, and 
$\{\omega_i^\alpha\}$ which are defined in the 
frame bundles of first, second, and third order, respectively. The 
forms $\omega_\alpha^i$ define a displacement of a point $x$ 
along the manifold $M$ on which the almost Grassmann structure 
is defined. By  (4.35), 
for $\omega_\alpha^i = 0$ these forms satisfy the equations 
\begin{equation}\label{eq:4.36}
d \pi_\alpha^\beta = \pi_\alpha^\gamma \wedge \pi_\gamma^\beta, 
\;\; d \omega_j^i = \omega_j^k \wedge \omega_k^i,  
\end{equation}
\begin{equation}\label{eq:4.37}
d \pi = 0, 
\end{equation}
\begin{equation}\label{eq:4.38}
d \pi_i^\alpha = \pi_i^\beta \wedge \omega_\beta^\alpha 
+  \pi_i^l \wedge  \pi_l^\alpha,
\end{equation}
and  the forms   $\pi_\beta^\alpha$ and $\pi_i^j$ satisfy the 
conditions (3.13), 
$$
\pi_\delta^\delta = 0,  \;\; \pi_l^l = 0.
$$

In view of (4.37), the form $\pi$ is an invariant form of 
the group ${\bf H}$ of homotheties of the tangent space 
$T_x (M)$, and the forms $\pi_j^i,  \pi_\beta^\alpha$, and 
$\pi$ are invariant forms 
of the group $G \cong {\bf SL} (p) \times {\bf SL} (q) \times 
{\bf H}$, whose transformations leave the cone 
$SC_x (p, q) \subset T_x (M)$ invariant. The group $G$ is the 
structural group of the almost Grassmann structure 
$AG (p-1, p+q-1)$. 

Equations (4.36), (4.37), and (4.38) prove that the forms 
$\pi, \pi_j^i,  \pi_\beta^\alpha$, and  $\pi_i^\alpha$ are 
invariant forms of the group $G'$ which is obtained as a 
differential prolongation of the group $G$. The group $G'$ is 
isomorphic to the group $G \htimes {\bf T} (pq)$ 
whose subgroup ${\bf T} (pq)$ is defined by the invariant forms  
$\pi_i^\alpha$.

To describe the group $G'$ geometrically, we compactify 
the tangent subspace $T_x (M)$ by enlarging it by the point 
at infinity and the Segre cone $SC_\infty (p, q)$ with its vertex 
at this point. Then the manifold $T_x (M) \cap SC_\infty (p, q)$ 
is equivalent to the algebraic variety $\Omega (p-1, p+q-1)$. 
Since the point $x$ at which the variety $\Omega (p-1, p+q-1)$ is 
tangent to the manifold $M$ is fixed, the geometry defined by the 
group $G'$ on  $\Omega (p-1, p+q-1)$ is equivalent to that of the 
flat Segre-affine space $SA^\rho$ of dimension $\rho = pq$, on 
which  the variety $\Omega (p-1, p+q-1)$  
is projected by means 
of a stereographic projection from the point $x$ (see 
\cite{SR85}, \cite{S31} or \cite{AG96}, \S  6.6). The group $G'$ 
is the group of motions of this space, its subgroup $G$ is the 
isotropy group of this space, and the subgroup ${\bf T} (pq)$ is 
the subgroup of parallel translations. We have already discussed 
this at the end of subsection {\bf 1.1}.

{\bf 4.} After the first reduction of second-order frames 
associated with the almost Grassmann structure 
$AG (p - 1, p + q - 1)$, we introduced the first structure tensor 
(the torsion tensor) $a = \{a^{i\beta\gamma}_{\alpha jk}\}$ of 
$AG (p - 1, p + q - 1)$. Let us set 
$b^1 = \{b_{ijk}^{l\beta\gamma}\},  
b^2 = \{b_{\delta jk}^{\alpha\beta\gamma}\}$, 
and $b = (b^1, b^2)$. 

Equations (4.27) and (4.28) show that the quantities 
$(a, b^1)$ and $(a, b^2)$ form linear homogeneous 
objects. They represent two subobjects of the {\em second 
structure object} $(a, b)$  of the  almost Grassmann 
structure $AG (p-1, p+q-1)$.

Taking the exterior derivatives of the  last equation 
of (4.35), we arrive at the following exterior cubic equation: 
\begin{equation}\label{eq:4.39}
\renewcommand{\arraystretch}{1.5}
\begin{array}{ll}
\bigl[\nabla c_{ijk}^{\alpha \beta \gamma}  
&\!\!\!\!\!\! + 3 c_{ijk}^{\alpha \beta \gamma}  \omega 
- b_{\sigma kj}^{\alpha \gamma \beta} \omega_i^\sigma 
+ b_{ikj}^{s \gamma \beta} \omega_s^\alpha 
+ (a^{s \alpha\beta}_{\sigma ij} 
a^{m\sigma\gamma}_{\varepsilon sk} 
+ a^{m \alpha\sigma}_{\varepsilon is} 
a^{s\gamma\beta}_{\sigma kj}) \omega^\varepsilon_m \\
&\!\!\!\!\!\! + (2 c_{ims}^{\alpha\varepsilon\sigma} 
a_{\sigma kj}^{s\gamma\beta}  \omega_\varepsilon^m 
- c_{sjk}^{\sigma\beta\gamma} 
a_{\sigma im}^{s \alpha\varepsilon}) 
\omega_\varepsilon^m\bigr] \wedge \omega^k_\gamma \wedge \omega^j_\beta = 0,
\end{array}
\renewcommand{\arraystretch}{1}
\end{equation}
where $\nabla c_{ijk}^{\alpha \beta \gamma} 
= d c_{ijk}^{\alpha \beta \gamma}  
-  c_{sjk}^{\alpha \beta \gamma}  \omega_i^s 
-  c_{isk}^{\alpha \beta \gamma}  \omega_j^s 
 -  c_{ijs}^{\alpha \beta\gamma}  \omega_k^s 
+  c_{ijk}^{\sigma \beta \gamma}  \omega_\sigma^\alpha 
+  c_{ijk}^{\alpha \sigma \gamma}  \omega_\sigma^\beta 
+  c_{ijk}^{\alpha \beta \sigma}  \omega_\sigma^\gamma$. 
For $\omega_\alpha^i = 0$, it follows  
from equation (4.39) that  
\begin{equation}\label{eq:4.40}
\nabla_\delta c_{ijk}^{\alpha \beta \gamma}  
+ 3 c_{ijk}^{\alpha \beta \gamma}  \pi 
- b_{\sigma kj}^{\alpha \gamma \beta} \pi_i^\sigma 
+ b_{ikj}^{s \gamma \beta} \pi_s^\alpha 
+(a_{\sigma ij}^{s\alpha\beta} a_{\varepsilon sk}^{m\sigma\gamma} 
+ a_{\varepsilon is}^{m\alpha\sigma} 
a_{\sigma kj}^{s\gamma\beta}) \pi_m^\varepsilon = 0.
\end{equation}
Let us set $c = \{c^{\alpha\beta\gamma}_{ijk}\}$. 
Equations (4.40), (4.27), (4.28), and  (3.32)  
prove that $S = (a, b, c)$ 
 form a linear homogeneous object, which is called the {\em 
third structure object} of the almost Grassmann structure 
$AG (p-1, p+q-1)$. It is defined in a fourth-order differential 
neighborhood of $AG (p-1, p+q-1)$. As we proved earlier, 
its subobject $a$ is a relative tensor (the torsion tensor) 
 defined in a second-order differential 
neighborhood of $AG (p-1, p+q-1)$, and the subobjects 
$(a, b^1), (a, b^2)$, and $(a, b)$ are defined in a third-order 
differential neighborhood of $AG (p-1, p+q-1)$. 

 The third structural object $S = (a, b, c)$ 
is the {\em complete geometric object} of the almost Grassmann 
structure $AG (p - 1, p + q - 1)$, since if we prolong 
the structure equations (4.35) of $AG (p - 1, p + q - 1)$, 
all newly arising objects are expressed in terms of 
the components of the object $S$ and their Pfaffian derivatives.  
This follows from Theorem 4.1.

{\bf 5.} Now we will find new closed form equations and 
differential equations that the components of  $b^1$ and $b^2$ 
satisfy. 
First, note that since the object $a$ is a relative tensor 
of weight $-1$, its components satisfy the equations 
\begin{equation}\label{eq:4.41}
\nabla a^{i\beta\gamma}_{\alpha jk} 
+ a^{i\beta\gamma}_{\alpha jk} \omega 
= a^{i\beta\gamma\delta}_{\alpha jkl} \omega_\delta^l. 
\end{equation}
This equation is equivalent to the equation (3.32). 
The quantities 
$a^{i\beta\gamma\delta}_{\alpha jkl}$ occurring in 
equations (4.41) are Pfaffian derivatives of the components of 
the tensor $a$ with respect to the basis forms $\omega_\delta^l$. 
With respect to the indices $i, \beta, \gamma, \alpha, j, k$, 
they satisfy the same relations (3.30)--(3.31) as 
the tensor $a$. 

 Substituting for $\nabla a^{i\beta\gamma}_{\alpha jk} 
+ a^{i\beta\gamma}_{\alpha jk} \omega$ and for the 
forms $\Omega_\alpha^\beta, 
\Omega_j^i$, and $d\omega$ in equations (3.14) 
 their values from expansions (4.41)  and  (4.35) 
and equating to zero the coefficients of 
$ \omega_{\delta}^{l} \wedge \omega_{\beta}^{j} \wedge 
 \omega_{\gamma}^{k}$, we obtain the following equations: 
\begin{equation}\label{eq:4.42}
\delta_{[j}^{i} b_{|\alpha|kl]}^{[\beta\gamma\delta]} 
- \delta_{\alpha}^{[\beta} b_{[jkl]}^{|i|\gamma\delta]} 
+ a_{\alpha [jkl]}^{i[\beta\gamma\delta]} 
+ 2 a_{\alpha [j|m}^{i [\beta|\varepsilon} 
a_{\epsilon|kl]}^{m|\gamma\delta]} =0.
\end{equation}
As earlier, in this formula, the alternation is carried 
over the pairs of indices $
{\beta \choose j}, 
{\gamma \choose k}$, and $
{\delta \choose l}$.
 Equation (4.42) can be written in the form 
\begin{equation}\label{eq:4.43}
\delta_{[j}^{i} b_{|\alpha|kl]}^{[\beta\gamma\delta]} 
- \delta_{\alpha}^{[\beta} b_{[jkl]}^{|i|\gamma\delta]} 
= A_{\alpha jkl}^{i\beta\gamma\delta}, 
\end{equation}
where the quantities $A_{\alpha jkl}^{i\beta\gamma\delta}$ are skew-symmetric with respect to the last three pairs of indices 
and are expressed in terms of the components of 
the tensor $a$ and their Pfaffian derivatives.

We will now prove 
the following result:

\begin{theorem} 
 For $p > 2$ and $q > 2$, 
 the components of $b^2$ and $b^1$ are 
expressed in terms of the components of the tensor $a$ and their 
Pfaffian derivatives. 
\end{theorem}

{\sf Proof}. In fact the components 
of $b^1$ and $b^2$ satisfy equations (4.43) which are a 
nonhomogeneous system of linear equations with respect to 
the quantities $b_{\alpha km}^{\beta \gamma \delta}$ and 
$b_{ijk}^{l\beta\gamma}$. Consider the homogeneous system 
corresponding to this  nonhomogeneous system; that is, set 
$a_{\alpha jk}^{i\beta \gamma} = 0$ in this nonhomogeneous 
system. This gives  
$$
\delta_{[j}^{i} b_{|\alpha|lm]}^{[\beta\delta\varepsilon]} 
- \delta_{\alpha}^{[\beta} b_{[jlm]}^{|i|\delta\varepsilon]} 
 = 0,
$$
or
$$
\delta_{j}^{i} b_{\alpha lm}^{\beta\delta\varepsilon} 
+ \delta_{l}^{i} b_{\alpha mj}^{\delta\varepsilon\beta} 
+ \delta_{m}^{i} b_{\alpha jl}^{\varepsilon\beta\delta} 
- \delta_{\alpha}^{\beta} b_{jlm}^{i \delta\varepsilon} 
- \delta_{\alpha}^{\delta} b_{lmj}^{i \varepsilon\beta} 
- \delta_{\alpha}^{\varepsilon} b_{mjl}^{i \beta \delta} 
 = 0.
$$
Contracting the latter equations with respect to 
the indices $i$ and $l$, 
$i$ and $j$, and $i$ and $m$ and applying conditions (4.31), 
we obtain
\begin{equation}\label{eq:4.44}
\renewcommand{\arraystretch}{1.5}
\left\{
\begin{array}{ll}
q b_{\alpha mj}^{\delta\varepsilon\beta} 
+ b_{\alpha jm}^{\beta\delta\varepsilon} 
+ b_{\alpha jm}^{\varepsilon\beta\delta} 
= 0, \\ 
 b_{\alpha mj}^{\delta\varepsilon\beta} 
+ q b_{\alpha jm}^{\beta\delta\varepsilon} 
+ b_{\alpha mj}^{\varepsilon\beta\delta} 
= 0, \\ 
b_{\alpha jm}^{\delta\varepsilon\beta} 
+  b_{\alpha jm}^{\beta\delta\varepsilon} 
+ q b_{\alpha mj}^{\varepsilon\beta\delta} 
= 0.  
\end{array}
\right.
\renewcommand{\arraystretch}{1}
\end{equation}

If we symmetrize equations (4.44) with respect to the indices 
$j$ and $m$, we obtain a homogeneous system of equations 
with respect to  $b_{\alpha (jm)}^{\delta\varepsilon\beta}, \; 
b_{\alpha (jm)}^{\beta\delta\varepsilon}$, and 
$b_{\alpha (jm)}^{\varepsilon\beta\delta}$. The determinant of 
the matrix of coefficients of this system is equal to 
$$
\left| 
\begin{array}{rrr} 
q &1 & 1 \\ 
1 &q & 1 \\ 
1 & 1 & q 
\end{array}
\right| 
= (q-1)^2 (q+2),
$$ 
and it does not vanish if $q \geq 2$. Hence this system has only 
the trivial solution: $b_{\alpha (mj)}^{\delta\varepsilon\beta} 
= 0$.

In the same manner, if we 
alternate equations (4.44) with respect 
to the indices $j$ and $m$, we obtain a homogeneous system of 
equations with respect to 
$b_{\alpha [jm]}^{\delta\varepsilon\beta}, \;\; 
b_{\alpha [jm]}^{\beta\delta\varepsilon}$, and  
$ b_{\alpha [jm]}^{\varepsilon\beta\delta}$. 
The determinant of the matrix of coefficients 
of this system is equal to 
$$
\left| 
\begin{array}{rrr} 
 -q &1 & 1 \\ 
  1 &q & 1 \\ 
  1 &1 & -q  
\end{array}
\right| 
= (q+1)^2 (q-2),
$$ 
and it does not vanish if $q > 2$. Hence this system has only 
the trivial solution: $b_{\alpha [jm]}^{\delta\varepsilon\beta} 
= 0$. As a result the homogeneous system in question has only 
the trivial solution: $b_{\alpha jm}^{\delta\varepsilon\beta} 
= 0$ provided that $q > 2$; thus the original nonhomogeneous 
system has a unique solution expressing the quantities 
$b_{\alpha km}^{\beta \gamma \delta}$  in terms of the 
components  
$a_{\alpha jk}^{i\beta \gamma}$ of the tensor $a$ and 
their Pfaffian derivatives. 

In a similar manner we can prove that if $p > 2$, then 
the quantities  $b_{jkm}^{i  \alpha \beta}$ are expressed in 
terms of the components  $a_{\alpha jk}^{i\beta \gamma}$ 
of the  tensor $a$ and their Pfaffian derivatives. 
Note that the condition $q > 2$ is required only for 
finding of $ b_{\alpha [jm]}^{\varepsilon\beta\delta} 
= b_{\alpha jm}^{\varepsilon(\beta\delta)}$  and 
the condition $p > 2$ for finding of \linebreak 
$b_{jkm}^{l [\alpha\beta]} = b_{j(km)}^{l \alpha\beta}$. 
\rule{3mm}{3mm}

Now we can see that the tensor $a$  
satisfies certain differential equations. 
These equations can be obtained if we substitute for the 
components  of $b^1$ and $b^2$ in equations (4.43)
 their values found in the way 
indicated above. 
The  conditions obtained in this manner 
are analogues of the Bianchi 
equations in the theory of spaces with affine connection.

{\bf 6.} 
 Next we will find new closed form equations and 
differential equations that the components of  $c$  satisfy. 
If we substitute for the 1-forms 
$\omega_{ij}^{\alpha\beta}$ in equations (4.9) and (4.10) 
their values taken from (4.32) and apply (4.34), we arrive 
at the following exterior cubic equations:
\begin{equation}\label{eq:4.45}
\renewcommand{\arraystretch}{1.5}
 \begin{array}{ll}
\Delta b^{\beta\gamma\delta}_{\alpha kl} 
\wedge \omega_\gamma^k  \wedge \omega^l_\delta 
- \Bigl(\displaystyle \frac{pq}{p+q}  
\delta_\alpha^\varepsilon \delta_m^s 
c_{skl}^{\beta\gamma\delta}
+ 2 b^{\beta\gamma\sigma}_{\alpha ks} 
a^{s \delta\varepsilon}_{\sigma lm}\Bigr)  
 \omega_\varepsilon^m
 \wedge \omega_\gamma^k  \wedge \omega_\delta^l = 0,\\
\end{array} 
\renewcommand{\arraystretch}{1}
\end{equation}
\begin{equation}\label{eq:4.46}
\renewcommand{\arraystretch}{1.5}
 \begin{array}{ll}
\Delta  b^{i\gamma\delta}_{jkl} 
\wedge \omega_\gamma^k  \wedge \omega^l_\delta 
+ \Bigl(\displaystyle \frac{pq}{p+q}  
\delta_m^i \delta^\varepsilon_\sigma  
c^{\sigma\gamma\delta}_{jkl} 
+ 2 b^{i\sigma\varepsilon}_{jsm} 
a^{s \gamma\delta}_{\sigma kl}\Bigr)  \omega_\varepsilon^m 
 \wedge \omega_\gamma^k  \wedge \omega_\delta^l = 0,
\end{array} 
\renewcommand{\arraystretch}{1}
\end{equation}
where 
$$
\renewcommand{\arraystretch}{1.5}
 \begin{array}{ll}
\Delta b^{\beta\gamma\delta}_{\alpha kl} 
= \nabla b^{\beta\gamma\delta}_{\alpha kl} 
+ 2 b^{\beta\gamma\delta}_{\alpha kl} \omega 
+ \displaystyle \frac{pq}{p+q} \Bigl(\delta_\alpha^\gamma 
\delta_k^s a_{\varepsilon sl}^{m\beta\delta} 
- \delta_\varepsilon^\beta \delta_s^m 
a_{\alpha kl}^{s\gamma\delta} \Bigr) \omega_m^\varepsilon, \\
 \Delta  b^{i\gamma\delta}_{jkl} 
= \nabla b^{i\gamma\delta}_{jkl}+ 2b^{i\gamma\delta}_{jkl} \omega 
+ \displaystyle \frac{pq}{p+q} \Bigl(\delta_\varepsilon^\sigma 
\delta_j^m a_{\sigma kl}^{i\gamma\delta} 
- \delta_\sigma^\gamma \delta_k^i 
a_{\varepsilon jl}^{m\sigma\delta} \Bigr)  \omega_m^\varepsilon. 
\end{array} 
\renewcommand{\arraystretch}{1}
$$
It follows from equations (4.45) and (4.46) that 
\begin{equation}\label{eq:4.47}
\Delta  b^{i\gamma\delta}_{jkl} 
= b^{i\gamma\delta\varepsilon}_{jklm} \omega_\varepsilon^m, \;\;
\Delta b^{\beta\gamma\delta}_{\alpha kl} 
= b^{\beta\gamma\delta\varepsilon}_{\alpha klm} 
\omega_\varepsilon^m, 
\end{equation}
where $b^{i\gamma\delta\varepsilon}_{jklm}$ and 
$b^{\beta\gamma\delta\varepsilon}_{\alpha klm}$ 
are the 
Pfaffian derivatives of $b^{\beta\gamma\delta}_{\alpha kl}$ 
and $b^{\beta\gamma\delta}_{\alpha kl}$, 
respectively. Substituting (4.47) into equations (4.45) and 
(4.46), we find the following differential equations for 
the components of $b$: 
\begin{equation}\label{eq:4.48}
\renewcommand{\arraystretch}{1.6}
 \begin{array}{ll}
 b^{\beta[\gamma\delta\varepsilon]}_{\alpha [klm]} 
- \displaystyle \frac{pq}{p+q}  \delta_\alpha^{[\varepsilon}
  c_{[mkl]}^{|\beta|\gamma\delta]}
- 2 b^{\beta[\gamma|\sigma}_{\alpha [k|s} 
a^{s| \delta\varepsilon]}_{\sigma| lm]} = 0,\\
 b^{i[\gamma\delta\varepsilon]}_{j[klm]} 
+ \displaystyle \frac{pq}{p+q}  
  \delta_{[m}^i c^{[\varepsilon\gamma\delta]}_{|j|kl]} 
+ 2 b^{i\sigma[\varepsilon}_{js[m} 
a^{|s| \gamma\delta]}_{|\sigma| kl]} = 0.
\end{array} 
\renewcommand{\arraystretch}{1}
\end{equation}

 Equations (4.48) can be written in the form 
\begin{equation}\label{eq:4.49}
  \delta_\alpha^{[\varepsilon}
  c_{[mkl]}^{|\beta|\gamma\delta]} 
= B^{\beta\gamma\delta\varepsilon}_{\alpha klm} 
\end{equation}
and 
\begin{equation}\label{eq:4.50}
  \delta_{[m}^i c^{[\varepsilon\gamma\delta]}_{|j|kl]} 
=  B^{i\gamma\delta\varepsilon}_{jklm}, 
\end{equation}
where the quantities 
$B^{\beta\gamma\delta\varepsilon}_{\alpha klm}$ and 
$B^{i\gamma\delta\varepsilon}_{jklm}$ 
are skew-symmetric with respect to the last three pairs of 
indices 
and are expressed in terms of the components of 
the subobjects $(a, b^1)$ and $(a, b^2)$, respectively, 
and  their Pfaffian derivatives.

We will now prove the following result:

\begin{theorem}  
If $p > 2$, then 
the components of $c$  are 
expressed in terms of the components of the subobject $(a, b^1)$ 
and their Pfaffian derivatives, and  if $q > 2$, then 
the components of $c$  are 
expressed in terms of the components of the subobject $(a, b^2)$ 
and their Pfaffian derivatives. 
\end{theorem} 

{\sf Proof.} 
We will prove only the first part of this theorem. The proof of 
the second part is similar. 

The components 
of $c$  satisfy equations (4.49) that are a 
nonhomogeneous system of linear equations with respect to 
$c_{jkl}^{\beta \gamma \delta}$. 
Consider the homogeneous system 
corresponding to this  nonhomogeneous system; that is, set 
$a = b^1 = 0$ in this nonhomogeneous 
system. This gives  
$$
\delta_\alpha^\varepsilon c_{jkl}^{\beta\gamma\delta} 
+ \delta_\alpha^\delta c_{ljk}^{\beta\varepsilon\gamma} 
+ \delta_\alpha^\gamma c_{jlk}^{\beta\delta\varepsilon}  = 0.
$$
Contracting this equation with respect to the indices 
$\alpha$ and $\varepsilon$, $\alpha$, and $\delta$, and 
$\alpha$ and $\gamma$, we obtain
\begin{equation}\label{eq:4.51}
\renewcommand{\arraystretch}{1.5}
\left\{
\begin{array}{ll}
p c_{jkl}^{\beta\gamma\delta} 
+ c_{ljk}^{\beta\delta\gamma} 
+ c_{klj}^{\beta\delta\gamma} = 0, \\
 c_{jkl}^{\beta\gamma\delta} 
+ p c_{ljk}^{\beta\delta\gamma} 
+ c_{klj}^{\beta\gamma\delta} = 0, \\
 c_{jkl}^{\beta\gamma\delta} 
+ c_{ljk}^{\beta\gamma\delta} 
+ p c_{klj}^{\beta\delta\gamma} = 0. 
\end{array}
\right.
\renewcommand{\arraystretch}{1}
\end{equation}

If we symmetrize and alternate 
equations (4.51) with respect to the indices 
$\gamma$ and $\delta$, we obtain two 
homogeneous systems of equations 
with respect to  $c_{jkl}^{\beta(\gamma\delta)} $ and 
$c_{jkl}^{\beta[\gamma\delta]}$ with different order 
of lower indices (cf. Subsection {\bf 4.5}). 
The determinants of 
the matrices of coefficients of these systems are equal to 
 $(p-1)^2 (p+2)$ and $(p+1)^2 (p -2)$, respectively. 
They do not vanish if $p > 2$. Hence these systems have only 
the trivial solution.  As a result the homogeneous system in 
question has only 
the trivial solution $c_{jkl}^{\beta\gamma\delta} = 0$ 
provided that $p > 2$; thus the original nonhomogeneous 
system has a unique solution expressing the components of  
$c$  in terms of the components of the subobject 
$(a, b^1)$ and their Pfaffian derivatives. \rule{3mm}{3mm} 

Now we can see that the object $(a, b)$  
satisfies certain differential equations. 
These equations can be obtained if we substitute for the 
components  of $c$  in equations (4.49) and (4.50) 
 their values found in the way indicated above. 
The  conditions obtained are other analogues of the Bianchi 
equations in the theory of spaces with affine connection.

{\bf 7.} An almost  Grassmann structure $AG (p-1, p+q-1)$ 
is said to be {\em locally Grassmann} (or {\em locally flat}) 
if it is locally equivalent to a Grassmann structure. 
This means that a locally flat 
almost  Grassmann structure $AG (p-1, p+q-1)$ admits 
a mapping onto an open domain of the algebraic variety 
$\Omega (m, n)$ of a projective space 
$P^N$, where $N = 
{n+1 \choose m+1} - \nolinebreak 1, \linebreak 
m = p - 1, n = p + q - 1$, under which the Segre cones of 
the  structure $AG (p-1, p+q-1)$ correspond to the asymptotic 
cones of the variety $\Omega (m, n)$.

 From the equivalence theorem of \'{E}. Cartan 
(see \cite{C08} or  \cite{Ga89}), it follows that in 
order for an almost  Grassmann structure $AG (p-1, p+q-1)$ 
to be locally Grassmann, it is necessary and sufficient that 
its structure equations have the form (1.4), (1.10), 
(1.11), and (1.12). Comparing these equations with equations 
(4.35), we see that an almost  Grassmann structure 
$AG (p-1, p+q-1)$ is locally Grassmann if and only if its 
complete structure object $S = (a, b, c)$ vanishes.

However, we established in this section that if $p > 2$ and 
$q > 2$, the components of $b$ are expressed in terms of the 
components of the tensor $a$ and their Pfaffian derivatives, 
and the components of $c$ are expressed in terms of the 
components of the subobject $(a, b)$ and their Pfaffian 
derivatives. Moreover it follows from our considerations 
that the vanishing of the tensor $a$ on 
a manifold $M$ carrying an almost  Grassmann structure 
implies the vanishing of the components of $b$ and $c$. 

Thus we have proved the following result: 

\begin{theorem} For $p > 2$ and $q > 2$,  an almost  Grassmann structure  $AG (p-1, p+q-1)$ is locally Grassmann if and only if 
its first structure tensor $a$ vanishes. 
\end{theorem}

Note that in 
the main parts of \cite{D93} and \cite{D94} the author 
considered torsion-free almost Grassmann  structures. He    
did not have results of our Theorems {\bf 4.2, 4.3} 
and {\bf 4.4} and erroneously assumed that for 
$p > 2, q > 2$ there exist torsion-free almost Grassmann  
structures that are not locally Grassmann  (locally flat) 
structures. According to Theorem {\bf 4.4}, if $p > 2$ and 
$q > 2$, then 
torsion-free almost Grassmann  structures are locally Grassmann  
structures. This is the reason that the results 
of  the main parts of \cite{D93} and \cite{D94} are valid 
only for $p = 2, q > 2; p > 2, q = 2$, and $p = q = 2$.   

{\bf 8.} We will now write  the structure equations 
(4.35) in index-free notation. To this end, we denote 
the matrix 1-form $(\omega_\alpha^i)$, defined in 
a first-order frame bundle, by $\theta$ and write 
equation (3.33) in the form 

\begin{equation}\label{eq:4.52}
 d\theta =  \omega \wedge \theta -  \theta  \wedge \omega_\alpha 
-  \omega_\beta \wedge \omega + \Theta, 
\end{equation}
where $\omega_\alpha = (\omega_\beta^\alpha)$ and 
$\omega_\beta = (\omega_j^i)$ are the matrix 1-forms 
defined in a second-order frame bundle for which 
$$
\mbox{{\rm tr}} \, \omega_\alpha = 0,  
\mbox{{ \rm tr}} \, \omega_\beta = 0;
$$
 the  form $\omega$ the scalar form  occurring in equation (3.33) 
and also defined in a second-order frame bundle.

The 2-form $\Theta = (\Theta_\alpha^i)$ is the {\em torsion form} 
with the components 
\begin{equation}\label{eq:4.53}
 \Theta_\alpha^i =  a^{i\beta\gamma}_{\alpha jk} \omega_\gamma^k 
 \wedge \omega_\delta^l. 
\end{equation}
The remaining structure equations (4.35) can be written 
in the form 
\begin{equation}\label{eq:4.54}
\renewcommand{\arraystretch}{1.5}
\begin{array}{ll}
d \omega_\alpha  + \omega_\alpha \wedge \omega_\alpha  
= \displaystyle \frac{q}{p+q} \Bigl[- I_\alpha \,
\mbox{{\rm tr}}\, (\varphi \wedge \theta) 
+ p \varphi \wedge \theta\Bigr] 
+ \Omega_\alpha,  \\
d \omega_\beta  + \omega_\beta \wedge \omega_\beta  
= \displaystyle \frac{p}{p+q} \Bigl[- I_\beta \,
\mbox{{\rm tr}}\, (\varphi \wedge \theta) 
+ q \theta \wedge \varphi\Bigr] 
+ \Omega_\beta,  \\
d \omega  =  \mbox{{\rm tr}}\, (\varphi \wedge \theta), \\
 d \varphi  + \omega_\alpha \wedge \varphi  + \varphi \wedge 
\omega_\beta  + \omega \wedge  \varphi = - (a \varphi) 
\wedge \theta + \Phi,
\end{array}
\renewcommand{\arraystretch}{1}
\end{equation}
where $\varphi = (\omega_i^\alpha)$ is a matrix 1-form 
defined in a third-order frame bundle;
 $I_\alpha = (\delta_\alpha^\beta)$ and $I_\beta = (\delta_i^j)$ are the unit tensors of orders $p$ and $q$, 
respectively; and 2-forms 
  $\Omega_\alpha = (\Omega^\alpha_\beta), \;  \Omega_\beta 
= (\Omega^i_j)$, and $\Phi = (\Phi^\alpha_i)$ are the 
 {\em curvature $2$-forms} of the 
 $AG (p-1, p+q-1)$-structure whose components are 
\begin{equation}\label{eq:4.55}
\Omega^\alpha_\beta = b_{\beta kl}^{\alpha \gamma \delta} \omega_\gamma^k \wedge \omega_\delta^l, \;\; 
\Omega^i_j = b_{jkl}^{i \gamma \delta} \omega_\gamma^k \wedge 
\omega_\delta^l, \;\; 
\Phi_i^\alpha = c_{ikl}^{\alpha \gamma \delta} 
\omega_\delta^l \wedge \omega_\gamma^k.
\end{equation}

\section{Semiintegrability of Almost Grassmann \newline Structures}

\setcounter{equation}{0}

{\bf 1.} In this section we will formulate 
 geometric conditions for an almost Grassmann structure 
$AG(p-1, p+q-1)$ defined on a manifold 
$M$ to be semiintegrable. The 
conditions will be expressed in terms of the complete structure object 
$S$ of  the almost Grassmann structure 
$AG (p-1, p+q-1)$ and its subobjects $S_\alpha$ and $S_\beta$ 
which will be defined in this section. For the proof of 
all theorems of this section see \cite{AG97} or  \cite{AG96}, 
\S 7.4.

First we will state  the theorem  
on the decomposition of the torsion tensor of an almost Grassmann 
structure $AG (p-1, p+q-1)$: 

\begin{theorem}
The torsion tensor $a  \! 
= \! \{a_{\alpha jk}^{i\beta\gamma}\}$ 
of the almost Grassmann structure 
$AG (p-1, p+q-1)$  decomposes into two subtensors:
\begin{equation}\label{eq:5.1}
a \!= \!a_\alpha \dot{+} a_\beta,
\end{equation}
where 
$$
a_\alpha =  \{a_{\alpha (jk)}^{i\beta\gamma}\}, \;\; 
a_\beta = \{ a_{\alpha jk}^{i(\beta\gamma)}\}.
$$
\end{theorem}

Note that the subtensors $a_\alpha$ and $a_\beta$ 
can be also represented in the form 
$$
a_\alpha =  \{a_{\alpha jk}^{i[\beta\gamma]}\}, \;\; 
a_\beta = \{ a_{\alpha [jk]}^{i\beta\gamma}\}.
$$

Note also that like the tensor $a$, its subtensors $a_\alpha$ and 
$a_\beta$ are skew-symmetric  with respect to 
the pairs of indices  $
{\beta \choose j}$ and 
$
{\gamma \choose k}$:
$$
a_{\alpha (jk)}^{i\beta\gamma} = -  
a_{\alpha (kj)}^{i\gamma\beta}, \;\; 
a_{\alpha [jk]}^{i\beta\gamma} = -  
a_{\alpha [kj]}^{i\gamma\beta},  
$$
  and they are also trace-free, 
since it follows from  (3.30) that 

\begin{equation}\label{eq:5.2}
a_{\alpha (jk)}^{i\alpha\gamma} = 0, \;\;
a_{\alpha ik}^{i[\beta\gamma]} = 0, \;\; 
a_{\alpha ik}^{i(\beta\gamma)} = 0, \;\;
a_{\alpha [jk]}^{i\alpha\gamma} = 0. 
\end{equation}

Next we will formulate the theorem on vanishing 
the subtensors $a_\alpha$ and $a_\beta$ for 
$p = 2$ and $q = 2$, respectively:

\begin{theorem}
If $p = 2$, then $a_\alpha  = 0$, 
and if $q = 2$, then $a_\beta = 0$. 
\end{theorem}

The following theorem gives  necessary and sufficient 
conditions for an 
almost Grassmann structure $AG (p-1, p+q-1)$ to 
be $\alpha$-semiintegrable or  $\beta$-semiintegrable.

\begin{theorem}

\begin{description}
\item[(i)] 
If $p > 2$ and $q \geq 2$, then for an almost Grassmann structure 
$AG (p-1, p+q-1)$ to be $\alpha$-semiintegrable, it is necessary 
and sufficient that the  condition $a_\alpha =
b_\alpha^1 = b_\alpha^2 = 0$ holds.

\item[(ii)] 
If $p \geq 2$ and $q > 2$, then for an almost Grassmann structure 
$AG (p-1, p+q-1)$ to be $\beta$-semiintegrable, it is necessary 
and sufficient that the condition $a_\beta =   
b_\beta^1 = b_\beta^2 = 0$ holds.

\end{description}
\end{theorem}

\vspace*{2mm}

We introduce the following notations:
$$
\begin{array}{lll}
b_\alpha^1 = \{b_{(jkl)}^{i\gamma\delta}\}, & \!\! 
b_\alpha^2 = \{b_{\alpha kl}^{[\beta\gamma\delta]}\}, & \!\! 
c_\alpha = \{c_{(ijk)}^{[\alpha\beta\gamma]}\}, \\ 
&&\\
b_\beta^1 = \{b_{[jkl]}^{i\gamma\delta}\}, & \!\! 
b_\beta^2 = \{b_{\alpha kl}^{(\beta\gamma\delta)}\}, & \!\! 
c_\beta = \{c_{[ijk]}^{(\alpha\beta\gamma)}\}. 
\end{array}
$$

For the cases when $p = 2$ or $q = 2$ or $p = q = 2$, 
the following theorem gives  
necessary and sufficient conditions for an 
almost Grassmann structure $AG (p-1, p+q-1)$ to 
be $\alpha$-semiintegrable or  $\beta$-semiintegrable: 

\begin{theorem} 
\begin{description}
\item[(i)] If $p = 2$, then the structure subobject $S_\alpha$ 
consists only of the tensor $b_\alpha^1$, and the vanishing of 
this tensor is necessary and 
sufficient for the almost Grassmann structure $AG (1, q + 1)$ 
to be $\alpha$-semiintegrable. 

\item[(ii)] If $q = 2$, then the structure subobject $S_\beta$ 
consists only of the tensor $b_\beta^2$, and the vanishing 
of this tensor is necessary and 
sufficient for the almost Grassmann structure $AG (p-1, p + 1)$ 
$($which is equivalent to the structure $AG (1, p + 1))$ 
to be $\beta$-semiintegrable. 

\item[(iii)] If $p = q = 2$, then the complete structural 
object $S$ consists only of the tensors 
$b_\alpha^1$ and $b_\beta^2$, and the vanishing 
of one of these tensors is necessary and 
sufficient for the almost Grassmann structure $AG (1, 3)$ 
to be $\alpha$- or $\beta$-semiintegrable, respectively. 
\end{description}
\end{theorem}

We conclude by the following remarks:

\begin{description}
\item[1.] The tensors 
$b_\alpha^1$ and $b_\beta^2$ are defined in a third-order 
differential neighborhood of the  almost Grassmann structure.

\item[2.] For $p = q = 2$, as was indicated earlier 
(see Example {\bf 2.3}), 
the almost Grassmann structure $AG (1, 3)$ is equivalent to the 
conformal $CO (2, 2)$-structure. Thus it follows 
(see \cite{AG96}, \S 5.1) 
that we have the following decomposition of its 
complete structural object: $S = b_\alpha^1 \dot{+} b_\beta^2$. 
This matches the splitting of the tensor of conformal curvature 
of the $CO (2, 2)$-structure: $C = C_\alpha \dot{+} C_\beta$. 
\end {description}

\end{document}